\documentclass[]{IEEEtran}
\usepackage{cite}

\usepackage{amsmath} 
\usepackage{amssymb}  

\usepackage{graphicx}

\usepackage{tikz}
\usepackage{pgfplots}
\pgfplotsset{compat=1.3}
\usetikzlibrary{shapes,arrows,arrows.meta}
\usepackage{arydshln} 

\usepackage{amsthm} 
\newtheorem{theorem}{Theorem}
\newtheorem{lemma}[theorem]{Lemma}
\newtheorem{corollary}[theorem]{Corollary}

\newtheorem{assumption}{Assumption}
\newtheorem{definition}{Definition}

\newtheorem{remark}{Remark}
\newtheorem{Problem}{Problem}

\usepackage[]{todonotes}
\usepackage{comment}

\usepackage{algorithmic}
\usepackage{graphicx}
\usepackage{textcomp}
\def\BibTeX{{\rm B\kern-.05em{\sc i\kern-.025em b}\kern-.08em
    T\kern-.1667em\lower.7ex\hbox{E}\kern-.125emX}}
\usepackage{fancyhdr,lipsum}
\makeatletter
\fancypagestyle{firstpage}{
	  \fancyhf{}
	  \fancyhead[C]{This is an author's copy intended to be submitted at IEEE and will be copyrighted by IEEE if accepted.}
	}
\begin{document}
\title{Robust Performance Analysis of Cooperative Control Dynamics via Integral Quadratic Constraints}
\author{Adwait~Datar, Christian~Hespe, and Herbert~Werner
\thanks{This work was supported by the German Research Foundation (DFG) within their priority programme SPP 1914 Cyber-Physical Networking.}
\thanks{Adwait Datar, Christian Hespe,  and Herbert Werner are with the Hamburg University of Technology, Institute of Control Systems, 21073 Hamburg, Germany. (e-mdail: \{adwait.datar, christian.hespe, h.werner\}@tuhh.de)}}

\maketitle
\thispagestyle{firstpage}
\begin{abstract}
We study cooperative control dynamics with gradient based forcing terms.
As a specific example, we focus on source-seeking dynamics with vehicles embedded in an unknown scalar field with a subset of agents having gradient information. As interaction mechanisms, formation control dynamics and flocking dynamics are considered.
We leverage the framework of $\alpha$-integral quadratic constraints to obtain convergence rate estimates whenever exponential stability can be achieved.
The communication graph and the interaction potential are assumed to be time-invariant and uncertain.
Sufficient conditions take the form of linear matrix inequalities independent of the size of network.
A derivation (purely in time-domain) of the so-called \textit{hard} Zames-Falb $\alpha$-IQCs involving general non-causal higher order multipliers is given along with a suitably adapted parameterization of the multipliers to the $\alpha$-IQC setting.
The time-domain arguments facilitate a straightforward extension to linear parameter varying systems.
Numerical examples illustrate the application of the theoretical results.
\end{abstract}

\begin{IEEEkeywords}
robust control, cooperative control, LMIs, 
\end{IEEEkeywords}
\section{Introduction}\label{sec:Introduction}
A variety of cooperative control algorithms (such as formation control, flocking, consensus, coverage control, etc.) lead to closed-loop dynamics which can be represented as gradient flows or dynamics involving gradient based forcing terms (see \cite{J.Cortes.2017}).
Extremum-seeking control problems \cite{michalowsky2016extremum} or problems containing real-time optimization \cite{nelson2018integral} lead to dynamics involving some gradient based forcing terms. 
The central theme of this paper is to analyze the performance of such dynamics using tools from robust control.
As a particular example, we focus on the source-seeking problem.
The abstract problem involves one or more vehicles located at arbitrary locations in an underlying scalar field with the goal of moving towards the minimum of the field which is called the source.
The gradient of the field at their respective locations is estimated in practice (\cite{khong2014multi, ogren2004cooperative, Datar.51220205152020}) and although the theoretical results are under the assumption of strictly or strongly convex fields and availability of gradients, experimental results \cite{Datar.51220205152020} with noisy measurements show that with a pre-filtering step, which involves fitting a strongly convex field with data collected in the local neighborhood, these protocols perform well in practice.
Our main tool for the analysis are integral quadratic constraints (IQCs) \cite{megretski1997system}
(See \cite{veenman2016robust, UlfJonsson.} for a tutorial).
\subsection{Related Work and Contributions}\label{sec:related_work} 
Stability analysis for source-seeking dynamics within the flocking framework is studied in \cite{Datar.51220205152020, Attallah.2020} using dissipativity based arguments by constructing valid (non-increasing) storage functions from physical energy-like functions.
 These physically motivated storage functions have a diagonal structure (such as the one used in \cite{OlfatiSaber.2006, Datar.51220205152020, Attallah.2020}).
The idea is to automate the search for storage functions (not necessarily diagonal) and include a less conservative stability analysis certificate along with performance guarantees.
These dissipativity based results can be interpreted as IQC results and the latter is the choice for the theoretical development in this paper.
The exponential versions of IQCs for systems in discrete-time ($\rho-$IQCs) are introduced in 
\cite{Lessard.2016} and in continuous-time ($\alpha-$IQCs) in \cite{Hu.2016}.
\cite{Hu.2016} introduces the \textit{soft} Zames-Falb (ZF) $\alpha$-IQCs corresponding to causal multipliers.
While \cite{Zhang.25022019} extends \cite{Lessard.2016} to less conservative non-causal ZF multipliers in the discrete-time setting, \cite{Freeman.62720186292018} presents the extension in the continuous-time setting.
The theory developed in \cite{Freeman.62720186292018} is in a very general setting of Bochner spaces and covers Lemma \ref{theom:lemma_pq_ZF} in this paper.
Moreover, proof of \cite[Lemma 3]{Freeman.62720186292018} (Lemma \ref{theom:lemma_pq_ZF} here) is not available and we therefore present a self-contained proof making all arguments in time-domain building on ideas presented in a sidebar in \cite{Scherer.16052021}.
The present paper goes further by considering multiplier parameterization proposed in \cite[Section 5.8.3]{veenman2016robust} adapted to the $\alpha$-IQCs setting.
\cite{Zhang.25022019} compares different multiplier factorizations which include the discrete-time analogue of the parameterization in \cite[Section 5.8.3]{veenman2016robust}.
In the spirit of \cite{Zhang.25022019, Freeman.62720186292018} which present examples showing the benefit of non-causal multipliers in the discrete-time case, we present an example of a continuous-time system with an integrator and an example of an LPV system demonstrating the benefit of non-causal higher order multipliers over causal ones.
Finally, previous works (such as \cite{Lessard.2016, Hu.2016, Freeman.62720186292018}) on the exponential version of IQCs present results for linear time invariant (LTI) systems.
Since we consider the \textit{hard}-IQCs with purely time-domain arguments, we show that it is rather straightforward to extend these results to linear parameter varying (LPV) \cite{shamma1992linear} systems as is done in \cite{pfifer2015robustness} for the standard IQCs.
This opens the doors for considering non-linear vehicle models with quasi-LPV representations.
The framework developed in this paper can be applied to extremum-seeking control problems \cite{michalowsky2016extremum} or problems involving realtime-optimization \cite{nelson2018integral}. 
Although the setup in \cite{nelson2018integral} is similar to the one in this paper, the focus there is on designing an optimizer and obtaining conditions for optimality and stability with static $\alpha-$IQCs.
The analysis of distributed optimization algorithms \cite{sundararajan2017robust, sundararajan2020analysis} leads to similar closed-loop dynamics.
In contrast to these works, we allow for a subset of agents to have the gradient information, relax the requirement on connectivity of the graph and instead require a connectivity to at least one agent with gradient information. 
As a result, the decomposition approach used in \cite{sundararajan2020analysis} does not go through as it is.
More importantly, since we do not look at time-varying graphs, we can use dynamic multipliers instead of the static multipliers used in \cite{sundararajan2020analysis}.
The analysis of flocking dynamics considered in this paper is closely related to the analysis of optimization algorithms over non-convex functions with non-exponential convergence rates. 
\cite{hu2017dissipativity} considers non-strongly convex functions whereas \cite{Fazlyab.2018} considers non-convex functions to obtain non-exponential convergence rates.
In contrast to \cite{Fazlyab.2018}, we use a time-invariant storage function and make use of LaSalle's invariance principle to conclude asumptotic stability.
The nature of the results obtained in \cite{sundararajan2020analysis} are very much in the same spirit as the one considered here.
The main contributions of this paper are now summarized:
\begin{enumerate}
    \item Building on \cite{Hu.2016}, a time-domain derivation of general non-causal exponential Zames-Falb (ZF) IQCs is given in Theorem \ref{theom:theorem_pqw_filtered_ZF} along with a suitable adaptation of multiplier parameterization from \cite[Section 5.8.3]{veenman2016robust} in Theorem \ref{theom:theorem_LMI_ZF}.  The key supporting Lemma \ref{theom:lemma_pq_ZF} which is at the core of this derivation is proven in Appendix \ref{appendix:Lemma1}.
	The present paper is an extension of \cite{datar2021iqc_source_seek} which does not contain proof of Lemma \ref{theom:lemma_pq_ZF}.
	
	\item Source-seeking dynamics involving vehicle agents under formation control are analyzed. 
	The equilibria of the dynamics are characterized in Lemmata \ref{lemm:minimizers} and \ref{lemm:minimizers_quadratic_radially_symmetric}. 
	LMIs (independent of network size) are derived for estimating the exponential convergence rates in {\color{black}{Lemma \ref{lemm:decomposition}}} and {\color{black}{Theorem \ref{theom:main_perf_analysis_formation}}} and supporting results are given in {\color{black}{Lemma \ref{lemm:f_in_Sml} and \ref{lemm:m_in_Sml_exists}}}.
	Numerical examples suggest that estimates obtained here are tight and are given in Section \ref{sec:numerical_results_formation_control}.
	
	\item Theorem \ref{theom:flocking_analysis_local} and Corollary \ref{corr:global_flocking_analysis} analyze source-seeking dynamics with flocking agents to give a local sufficient condition for asymptotic stability in the form of an LMI independent of network size generalizing \cite{Datar.51220205152020, Attallah.2020}.
	
	\item Numerous examples demonstrating the application of the theoretical results are given.
\end{enumerate}

\section{Notation}\label{sec:notation}
Let $\mathbb{R}$ denote the set of real numbers.
Let $\mathbb{S}^{n}$ denote the set of symmetric matrices of size $n$.
For $X \in \mathbb{S}^n$, $X\succ(\succeq)\;0$ means that $X$ is symmetric positive definite (semi-definite), and $X\prec(\preceq)\;0$ denotes that $-X$ is symmetric positive definite (semi-definite).
We use $*$ to denote appropriate entries to make the matrix symmetric. 
Let $\mathbf{0}$ and $\mathbf{1}$ denote the vectors or matrices of all zeroes and ones of appropriate sizes, respectively. 
Let $I_d$ be the identity matrix of dimension $d$ and we remove the subscript $d$ when the dimension is clear from context. 
Let $\otimes$ represent the Kronecker product.
For some $0<m \leq L$, let $\mathcal{S}(m,L)$ denote the set of continuously differentiable functions $f:\mathbb{R}^d\rightarrow \mathbb{R}$, for some positive integer $d$, which are strongly convex with parameter $m$, and have Lipschitz gradients with parameter $L$, i.e., for all $y_1, y_2 \in \mathbb{R}^d$,
\begin{equation*}
	m||y_1-y_2||^2 \leq  (\nabla f(y_1)-\nabla f(y_2))^T(y_1-y_2) \leq L ||y_1-y_2||^2.
\end{equation*}
The set of vector valued functions which are square-integrable over $[0,T]$ for any finite $T<\infty$ is denoted by $\mathcal{L}_{2e}[0,\infty)$.
Let $\mathcal{L}_1(-\infty,\infty)$ denote the set of functions $h:\mathbb{R}\rightarrow \mathbb{R}$, such that $\int_{-\infty}^{\infty}|h(t)|dt < \infty$.
We use 
$\left[\begin{array}{c|c}
	\mathcal{A}     &  \mathcal{B}\\
	\hline
	\mathcal{C}     &  \mathcal{D}
\end{array}\right]$
to represent an LTI system with state-space realization given by matrices $\mathcal{A},\mathcal{B},\mathcal{C},\mathcal{D}$.
Let the convolution operation be denoted by $a(t)*b(t)=\int_0^t a(t-\tau)b(\tau)d\tau$.
Interconnections between agents are modeled by graph structures.
The number of agents is denoted by $N$ and let $\mathcal{G}=(\mathcal{V},\mathcal{E})$ be an undirected unweighted graph of order $N$ with the set of nodes denoted by $\mathcal{V}=\left\{1,2,\ldots N \right\}$ and the set of edges by $\mathcal{E}\subseteq \mathcal{V}\times \mathcal{V}$.
The set of informed agents denoted by $\mathcal{V}_l \subseteq \mathcal{V}$ is a subset of nodes with the gradient information.
We use $\mathcal{G}_{\textnormal{star}}^N=(\{1,\cdots,N\},\{(1,1),(1,2),\cdots,(1,N)\})$ to denote a star and $\mathcal{G}_{\textnormal{cycle}}^N=(\{1,\cdots,N\},\{(1,2),(2,3),\cdots,(N-1,N),(N,1)\})$ to denote a cycle.
The notation $\hat{X}$ denotes the matrix $I_N \otimes X$ and let $X_{(d)}$ denote the matrix $X \otimes I_d$ for any matrix $X$.
For an ordered set of vectors $(x_1,x_2,\hdots,x_N)$, let the vector formed by stacking these vectors be denoted by $x=\left[x_1^T \; x_2^T \; \hdots \; x_N^T \right]^T$. 
\section{Problem Setup}\label{sec:problem_setup}
Consider a source-seeking scenario where $N$ vehicle agents moving in $\mathbb{R}^d$ space (typically $d\in\{1,2,3\}$) are embedded in an external differentiable scalar field $\psi: \mathbb{R}^d \xrightarrow[]{}\mathbb{R}$ with the goal of cooperatively moving towards the minimizer (source) of $\psi$.
To this end, assume that a local tracking controller has been designed and the closed-loop dynamics of the $i^{\textnormal{th}}$ vehicle agent with reference position $q_i(t)$ and reference velocity $p_i(t)$ can be described for a given initial condition $x_i(0)$ by 
\begin{equation}\label{eq:vehicle_dyn_flock_chap}
	\begin{split}
		\Dot{x}_i(t)&=Ax_i(t) + \begin{bmatrix} B_q & B_p \end{bmatrix} \begin{bmatrix} q_i(t)\\p_i(t)\end{bmatrix}, \\
		y_i(t)&=Cx_i(t),
	\end{split}
\end{equation}
where $y_i(t)$ is the position output of agent $i$.
These closed-loop dynamics are augmented by the second-order dynamics
\begin{equation}\label{eq:vir_vehicle_dyn_flock_chap}
	\begin{split}
		\Dot{q}_i(t)&=p_i(t),\\
		\Dot{p}_i(t)&=-k_d \cdot p_i(t) -k_p\cdot u_i(t),
	\end{split}
\end{equation}
where $u_i(t) \in \mathbb{R}^d$ is the external force, $q_i(0)=Cx_i(0)$, $p_i(0)=0$.
This is depicted in Fig. \ref{fig:plant G_flock_chap}. 
	\begin{figure}[t!]
		\centering
		\tikzstyle{block} = [draw, rectangle, 
    minimum height=3em, minimum width=6em]
\tikzstyle{sum} = [draw,  circle, node distance=1cm]
\tikzstyle{input} = [coordinate]
\tikzstyle{output} = [coordinate]
\tikzstyle{pinstyle} = [pin edge={to-,thin,black}, pin distance = 1.2cm]

\begin{tikzpicture}[auto, node distance=2cm,>=latex']
\node [input, name=input] {};
\node [block, right of=input,pin={[pinstyle]above:$\begin{bmatrix}q_i(0)\\p_i(0)\end{bmatrix}$},node distance=2.6cm](Flockcontrol) {$\begin{array}{cc}
   \Dot{q}_i =& p_i \\
   \Dot{p}_i =&-k_d p_i-k_pu_i
\end{array}$};
\node [block, right of=Flockcontrol,pin={[pinstyle]above:$x_i(0)$},node distance=4.2cm,minimum height=2.5em,minimum width=3em](velocityloop) {$\left[\begin{array}{c|c}
A     & \begin{bmatrix}B_q & B_p\end{bmatrix} \\
\hline
C     &  \mathbf{0}
\end{array}\right]$};
\node [output, right of=velocityloop,node distance=2.2cm](output) {};
\node[text=red, above left= 5mm and -1cm of Flockcontrol] (veh) {$G$};
\draw[thick,red, dashed] (veh.east)-|([xshift=-5.8mm]output.west)|-([yshift=-3mm]Flockcontrol.south)-|([xshift=0.6cm]input.east)|-(veh.west);
\draw [draw,->] (input) -- node {$u_i$} (Flockcontrol);
\draw [draw,->] (Flockcontrol) -- node {$\begin{bmatrix}q_i\\p_i\end{bmatrix}$} (velocityloop);
\draw [draw,->] (velocityloop) -- node[name=outmid] {$y_i$} (output);
\end{tikzpicture}
		\caption{Local control architecture on agent $i$.}	
		\label{fig:plant G_flock_chap}	
	\end{figure}
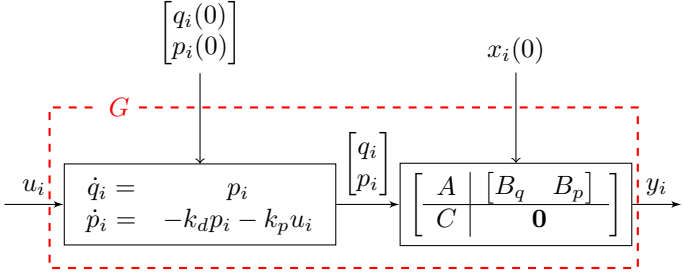
\begin{remark}\label{rem:pre-filter_motivation}
			The inclusion of pre-filter dynamics \eqref{eq:vir_vehicle_dyn_flock_chap} has multiple benefits. 
			First of all, these dynamics act as crude approximations of the vehicle dynamics generating smooth reference trajectories.
			If the vehicle is able to track the generated trajectories well enough, i.e., $y_i \approx q_i$, we can study the interconnected system of second order vehicle models and expect a stable behavior also with higher order vehicle models (see \cite{Datar.51220205152020, Attallah.2020}).
			This intuitive argument is made rigorous here.
			Furthermore, one can design tracking controllers by any method of choice and plug in the closed-loop dynamics in this framework for a systematic performance analysis.
			Finally, these pre-filter dynamics provide us with simple interpretable tuning knobs, e.g., the damping co-efficient $k_d$, which can be used to slow down or speed up the overall dynamics.
\end{remark}
The overall dynamics with state $\eta_i=\left[x_i^T \; q_i^T\; p_i^T\right]^T$ and suitable initial condition $\eta_i(0)$ can be represented by
\begin{equation} \label{eq:sys_dyn_G_flock_chap}
	\begin{split}
		\Dot{\eta}_i(t)&=A_G\eta_i(t) + B_G u_i(t), \\ 
		y_i(t)&=C_G \eta_i(t), \\
	\end{split}
\end{equation}
where 
\begin{equation*}
	\begin{split}
		\left[\begin{array}{c|c}
			A_G     &  B_G\\
			\hline
			C_G     &  D_G
		\end{array}\right]
		& =
		\left[\begin{array}{ccc|c}
			A            & B_q & B_p & \mathbf{0}\\
			\mathbf{0}   & \mathbf{0} & I_d & \mathbf{0}\\
			\mathbf{0}  & \mathbf{0} & -k_d I_d & -k_p I_d\\
			\hline
			C  & \mathbf{0} & \mathbf{0} & \mathbf{0}\\
		\end{array}
		\right].
	\end{split}
\end{equation*}

The interaction between different agents is captured by an interaction potential $V:\mathbb{R}^{Nd} \rightarrow \mathbb{R}$ (see \cite{OlfatiSaber.2006} or Section \ref{sec:numerical_results_flocking} for an example) and let a non-empty subset $\mathcal{V}_l \subseteq \mathcal{V}$ of informed agents have access to the local gradient $\nabla \psi$ evaluated at their respective positions.
These informed agents have an additional forcing term in the direction of the negative gradient that drives them towards the source.
For a compact representation of the overall dynamics, define a function $f:\mathbb{R}^{Nd}\rightarrow \mathbb{R}$ as follows.
\begin{definition}\label{defn:f_flock}
	For a given interaction potential $V$, the set of informed agents $\mathcal{V}_l$ and a scalar field $\psi$, define a function $f:\mathbb{R}^{Nd}\rightarrow \mathbb{R}$ by
	\begin{equation}\label{eq:defn_f_flock}
		f(y)= V (y)+ \sum_{i \in \mathcal{V}_l} \psi(y_i).
	\end{equation}
\end{definition}
\begin{remark}
			Note that $f$ is the sum of the interaction potential and the external scalar field and can be thought of as a sum of the cooperation cost and the source-seeking cost. 
\end{remark}

\begin{figure}[t!]
	\centering
	\def\svgwidth{0.5\textwidth}
	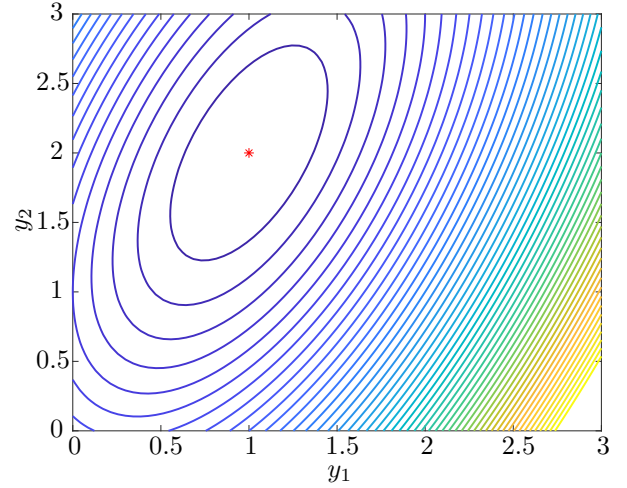	
	\caption{Contour plot of an example $f(y)= \frac{1}{2}(y_1+1-y_2)^2+ \frac{1}{2}(y_1-1)^2$ with $N=2$ and $d=1$, where minimizer $[1\,\,2]^T$ is denoted by red star.}	
	\label{fig:f_c_example}	
\end{figure}

\begin{figure}[t!]
	\centering
	\def\svgwidth{0.5\textwidth}
	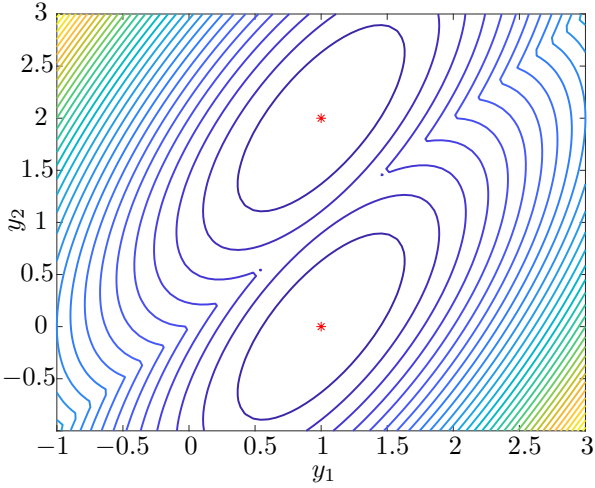
	\caption{Contour plot of an example $f(y)= \frac{1}{2}(|y_1-y_2|-1)^2+ \frac{1}{2}(y_1-1)^2$ with $N=2$ and $d=1$, where minimizers $[1\,\,0]^T$ and $[1\,\,2]^T$ are denoted by red stars.}	
	\label{fig:f_nc_example}	
\end{figure}

Figures \ref{fig:f_c_example} and \ref{fig:f_nc_example} show contour plots for two examples of $f$ corresponding to typical formation control and flocking interactions illustrating Definition \ref{defn:f_flock}. 
In both examples, we consider two connected agents moving in $\mathbb{R}$ with agent number 1 as the informed agent, i.e., $N=2$, $d=1$ and $\mathcal{V}_l=\{1\}$.
Furthermore, the external field $\psi(y_1)= \frac{1}{2}(y_1-1)^2$. 
For the first example corresponding to typical formation control dynamics, let the desired formation correspond to agent 2 being 1 unit ahead of agent 1 and this can be captured by a quadratic interaction $V(y)=\frac{1}{2}(y_1+1-y_2)^2$. 
We can write this using the graph Laplacian matrix $\mathcal{L}$ and formation reference vector $r=[0\quad 1]$ as 
\begin{equation*}
	\begin{split}
		V(y)&=\frac{1}{2}(y-r)^T(\mathcal{L} \otimes I_d)(y-r)\\
		&=\frac{1}{2}\begin{bmatrix}
			y_1&y_2-1
		\end{bmatrix}\begin{bmatrix}
			\,\,1&-1\\-1&\,\,1
		\end{bmatrix}\begin{bmatrix}
			y_1\\y_2-1
		\end{bmatrix}\\
		&=\frac{1}{2}(y_1+1-y_2)^2.
	\end{split}
\end{equation*}
Since agent $1$ is the informed agent, the desired arrangement corresponds to the minimum of $f$ at $y=[1\quad 2]^T$.
This is shown in Fig. \ref{fig:f_c_example}.
For the example corresponding to typical distance based interactions common in flocking, let the desired distance between the agents be 1 units which can be captured by $V(y)= \frac{1}{2}(|y_1-y_2|-1)^2$. 
Note that since agent $1$ is the informed agent, the desired arrangements which correspond to the minima of $f$ are $y=[1\quad 0]^T$ and $y=[1\quad 2]^T$, i.e., the informed agent at the source and agent $2$ at a distance of $1$ unit from the informed agent in either direction.
This is depicted in Fig. \ref{fig:f_nc_example}.

Using the notation introduced in Section \ref{sec:notation}, let $\eta(t), u(t)$ and $y(t)$ be obtained by stacking the states, inputs and outputs of the agents, respectively. The overall dynamics can be described as
\begin{equation} \label{eq:sys_dyn_G_hat}
	\begin{split}
		\Dot{\eta}(t)&=\hat{A}_G\eta(t) + \hat{B}_G u(t), \quad \quad \eta(0)=\eta_{0},\\
		y(t)&=\hat{C}_G \eta(t), \\
		u(t)&=\nabla f (y(t))=\nabla V(y(t)) + u_{\psi}(t),
	\end{split}
\end{equation}
where $u_{\psi}(t)$ is defined by stacking up
\begin{equation}\label{eq:forcing_term}
	\begin{split}
		u_{\psi_i}(t)&=
		\begin{cases}
			\nabla \psi(y_i(t)), & \text{if}\ i\in \mathcal{V}_l, \\
			0 & \text{otherwise}.
		\end{cases}
	\end{split}
\end{equation}
\subsection{Formation Control Dynamics with Convex Interactions}
	The special case of formation control dynamics \cite{Fax.2004} is obtained if $V(y)=\frac{1}{2}(y-r)^T(\mathcal{L} \otimes I_d)(y-r)$, where $\mathcal{L}$ is the graph Laplacian matrix and $r$ is the formation reference encoding the desired formation shape (see \cite{Fax.2004} for details).
	The objective is to make the agents form a desired shape and move towards the source of the external field at the same time.
	The assumptions and the first problem statement is given next.
	\begin{assumption}\label{assum:psi_Sml}
		The scalar field $\psi$ belongs to $\mathcal{S}(m_{\psi},L_{\psi})$ and let $y_{\textnormal{opt}}$ minimize $\psi$, i.e., $\psi(y)\geq \psi(y_{\textnormal{opt}}) \, \forall y \in \mathbb{R}^d$.
	\end{assumption}
	\begin{assumption}\label{assm:path_to_leaders}
		For every node $i \in \mathcal{V}$, there is a node $j \in \mathcal{V}_l$ such that $\mathcal{G}$ contains a path from $i$ to $j$.
	\end{assumption}
\begin{Problem} \label{prob:formation_control}
	Let $\psi$, $\mathcal{G}$ and $\mathcal{V}_l$ satisfy Assumptions \ref{assum:psi_Sml} and \ref{assm:path_to_leaders} and $V(y)=\frac{1}{2}(y-r)^T(\mathcal{L} \otimes I_d)(y-r)$. 
	Characterize the equilibria of dynamics \eqref{eq:sys_dyn_G_hat} and derive sufficient conditions independent of the network size $N$ under which the state trajectories generated by \eqref{eq:sys_dyn_G_hat} remain bounded for all $t\geq 0$ and $y$ converges exponentially with a rate $\alpha$ to the minimizer $y_*$ of $f$, i.e., $\exists \kappa \geq 0:\,||y(t)-y_*(t)||\leq \kappa e^{-\alpha t}\,\,\,\forall t\geq 0$.
\end{Problem}
\subsection{Flocking Dynamics with Non-Convex Interactions}
Flocking dynamics, such as from \cite{OlfatiSaber.2006}, involve distance based non-convex interaction potentials but satisfy a Lipschitz condition on the gradients. 
The objective here is to make the agents flock towards the source of the external field.
This requires allowing for a non-convex $f$ and prevents us from proving exponential stability with the help of ZF-IQCs (to be introduced in Section \ref{sec:IQC_theory}).
We therefore make less strict assumptions on $f$ (such as allowing non-convexity) and ask for asymptotic stability in the second problem statement.
\begin{assumption}\label{assum:f_flock}
	The function $f: \mathbb{R}^{Nd} \rightarrow \mathbb{R}$ and an open set $\mathcal{S}$ containing a local minimizer $y_*$ of $f$ satisfy the following conditions:
	\begin{itemize}
		\item[1)] $f$ is differentiable on $\mathcal{S}$. 
		\item[2)] $y_* \in \mathcal{S}$ is a local minimizer of $f$ on $\mathcal{S}$, i.e., $\nabla f(y_*)=0$ and $f(y)\geq f(y_*)=f_{\textnormal{min}}\quad \forall y  \in \mathcal{S}.$
		\item[3)] There exists a symmetric matrix $M_{10}$ such that
		\begin{equation*}
			(*)
			(M_{10} \otimes I_{Nd})
			\begin{bmatrix}
				y-y_*\\\nabla f(y)
			\end{bmatrix}
			\geq 0 \quad \forall y\in \mathcal{S}. 
		\end{equation*}
		
		\item[4)]There exists a symmetric matrix $M_{20}$ such that
		\begin{equation*}
			(*)
			(M_{20} \otimes I_{Nd})
			\begin{bmatrix}
				x-y\\\nabla f(x)-\nabla f(y)
			\end{bmatrix}
			\geq 0 \quad \forall x,y\in \mathcal{S}.
		\end{equation*}
		
		\item[5)] There exist constants $c_1$ and $c_2$ such that $$\mathcal{S}_0=\{y+e|\,f(y)-f_{\textnormal{min}}\leq c_1,||e||^2\leq c_2\}$$ is bounded and contained in $\mathcal{S}$.
	\end{itemize}
\end{assumption}
\begin{remark}
	The knowledge about $f$ needs to be captured via a suitable selection of matrices $M_{10}$ and $M_{20}$ in Assumption \ref{assum:f_flock} items 3) and 4). 
	See \cite[Section 6.2]{Fazlyab.2018} for a few examples.
	If global properties of $f$ are known, one can set $\mathcal{S}=\mathbb{R}^{Nd}$.
	Assumption \ref{assum:f_flock} item 5) is required to ensure an invariance property for the set $\mathcal{S}$ such that trajectories do not exit this set to regions of $\mathbb{R}^{Nd}$ where properties of $f$ are not known.
\end{remark}

An additional assumption on the vehicle dynamics is made which states that the local tracking controller is stabilizing and has zero steady state error for step position references $q(t)$.
\begin{assumption}\label{assm:ss_error}
	The vehicle dynamics \eqref{eq:vehicle_dyn_flock_chap} are such that every eigenvalue of the matrix $A$ has strictly negative real part and $-CA^{-1}B_q=I.$
\end{assumption}
\begin{Problem}\label{prob:flocking}
Let $f$ and the vehicle dynamics satisfy Assumptions \ref{assum:f_flock} and \ref{assm:ss_error}. Derive sufficient conditions independent of the network size $N$ under which the state trajectories generated by the dynamics \eqref{eq:sys_dyn_G_hat} remain bounded for all $t\geq 0$ and $y$ converges asymptotically to a minimizer of $f$.
\end{Problem}

\section{Analysis: Formation Control (Problem 1)} \label{sec:analysis_formation_control}
This section addresses Problem \ref{prob:formation_control}. 
Conditions ensuring $f \in \mathcal{S}(m,L)$ and equilibria of the dynamics (minimizers of $f$) are studied in Sections \ref{subsec:convex_smooth_f} and \ref{subsec:minimizers} followed by the main analysis in Section \ref{sec:IQC_theory} where $f \in \mathcal{S}(m,L)$ is assumed.
\subsection{Smoothness and Convexity Properties of $f$}\label{subsec:convex_smooth_f}
For convenience, let us define two grounded Laplacians.
\begin{definition}\label{defn:Lgs_defn}
	For a given graph $\mathcal{G}$ of order $N$ (with its corresponding Laplacian $\mathcal{L}$), a set of informed agents $\mathcal{V}_l$ and constants $0<m_{\psi}\leq L_{\psi}$, let the grounded Laplacians $\mathcal{L}_s$ and $\mathcal{L}_b$ be defined by
	\begin{equation} \label{eq:Lgs_defn}
		\begin{split}
			\mathcal{L}_s=\mathcal{L}+m_{\psi} E \textnormal{ and }       \mathcal{L}_b=\mathcal{L}+L_{\psi} E, \\
		\end{split}
	\end{equation}
	where $E$ is a diagonal matrix with the $i^{th}$ diagonal entry equal to $1$ if $i \in \mathcal{V}_l$ and equal to $0$ otherwise.
\end{definition}
\begin{remark}
			To compare the above with the definition of a \textit{grounded Laplacian} from literature, consider an undirected graph $\mathcal{G}_s$ by adding an $(N+1)^{\textnormal{th}}$ node to the vertex set of $\mathcal{G}$  such that this added node is grounded \cite{Xia.2017} and has an edge with all informed agents in $\mathcal{V}_l$ with an edge weight $m_{\psi}$. 
			The grounded Laplacian for $\mathcal{G}_s$ as defined in \cite{Xia.2017} equals $\mathcal{L}_s$.
\end{remark}
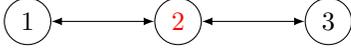
\begin{figure}[t!]
	\centering
	\begin{tikzpicture}	
		\node[draw, circle](1) at (0,0) {$1$};
		\node[draw, circle](2) at (2cm,0) {\textcolor{red}{$2$}};
		\node[draw, circle](3) at (4cm,0) {$3$};		
		\draw [draw,<->,>=latex] (1) -- (2);
		\draw [draw,<->,>=latex] (2) -- (3);
	\end{tikzpicture}
	\caption{Example graph with agent 2 as the informed agent for illustration of Definition \ref{defn:Lgs_defn} and Lemma \ref{lemm:f_in_Sml}.}
	\label{fig:Example_grounded_Lap}	
\end{figure}
To illustrate the above definition, consider a graph depicted in Fig. \ref{fig:Example_grounded_Lap} where the informed agent $2$ is shown in red. Furthermore, let constants $m_{\psi}=1$ and $L_{\psi}=2$.
The grounded Laplacians for this setup as defined above are given by
\begin{equation*}
	\begin{split}
		\mathcal{L}_s&=\mathcal{L}+1E=
		\begin{bmatrix}
			1 & -1 & 0\\
			-1 & 2+1 & -1\\
			0 & -1 & 1
		\end{bmatrix}=
	\begin{bmatrix}
	1 & -1 & 0\\
	-1 & 3 & -1\\
	0 & -1 & 1
\end{bmatrix}, \\ 
		\mathcal{L}_b&=\mathcal{L}+2E=		
		\begin{bmatrix}
			1 & -1 & 0\\
			-1 & 2+2 & -1\\
			0 & -1 & 1
		\end{bmatrix}==		
	\begin{bmatrix}
	1 & -1 & 0\\
	-1 & 4 & -1\\
	0 & -1 & 1
\end{bmatrix}.
	\end{split}
\end{equation*}
\begin{lemma} \label{lemm:f_in_Sml}
	For a given graph $\mathcal{G}$ of order $N$ (with its corresponding Laplacian $\mathcal{L}$), a set of informed agents $\mathcal{V}_l$, a scalar field $\psi$, a formation reference vector $r$ and constants $m_{\psi}$, $L_{\psi}$ such that $0<m_{\psi}\leq L_{\psi}$, let $f$ be as defined in Definition \ref{defn:f_flock} with $V(y)=\frac{1}{2}(y-r)^T(\mathcal{L} \otimes I_d)(y-r)$ and $\mathcal{L}_s$, $\mathcal{L}_b$ be as defined in Definition \ref{defn:Lgs_defn}.
	Then, for constants $m$, $L$ such that $0<m \leq L$, the following statements are equivalent:
	\begin{enumerate}
		\item[1)] $f\in \mathcal{S}(m,L) \textnormal{ for all } \psi \in \mathcal{S}(m_{\psi},L_{\psi})$,
		\item[2)] $m I \preceq  \mathcal{L}_s$ and $ \mathcal{L}_b \preceq L I$.
	\end{enumerate}
\end{lemma}
\begin{proof}
	See Appendix \ref{appendix:proofs}.
\end{proof}
The above lemma gives a condition to verify $f\in \mathcal{S}(m,L)$.
The direction $1) \implies 2)$ can be seen by picking a quadratic $\psi$ whereas the reverse direction is seen as a direct consequence of the definition of $\mathcal{S}(m,L)$.

Continuing our example from Fig. \ref{fig:Example_grounded_Lap} with $m_{\psi}=1$ and $L_{\psi}=2$, observe that $0.2679 \cdot I\preceq\mathcal{L}_s\preceq \mathcal{L}_b\preceq 4.5616 \cdot I$, which implies because of Lemma \ref{lemm:f_in_Sml} that $f\in \mathcal{S}(0.2679,4.5616) \textnormal{ for all } \psi \in \mathcal{S}(1,2)$.
The next lemma shows that Assumption \ref{assm:path_to_leaders} is necessary and sufficient for $f \in \mathcal{S}(m,L)$.
The main idea behind the proof is to see that $\mathcal{L}_s\succ 0$ if Assumption \ref{assm:path_to_leaders} is satisfied while it has a non-empty kernel if the assumption is not satisfied.
\begin{lemma}\label{lemm:m_in_Sml_exists}
	Let $0<m_{\psi}\leq L_{\psi}$ be fixed.
	There exist constants $m$, $L$ with $0<m\leq L$ such that $f$ as defined in Definition \ref{defn:f_flock} with $V(y)=\frac{1}{2}(y-r)^T(\mathcal{L} \otimes I_d)(y-r)$ belongs to $\mathcal{S}(m,L)$ for every $\psi \in \mathcal{S}(m_{\psi},L_{\psi})$ if and only if the graph $\mathcal{G}$ and the set of informed agents $\mathcal{V}_l$ satisfy Assumption \ref{assm:path_to_leaders}.
\end{lemma}
\begin{proof}
See Appendix \ref{appendix:proofs}.
\end{proof}
\subsection{Minimizers of $f$}\label{subsec:minimizers}

\begin{figure}[t!]
	\centering
	\includegraphics[scale=0.6]{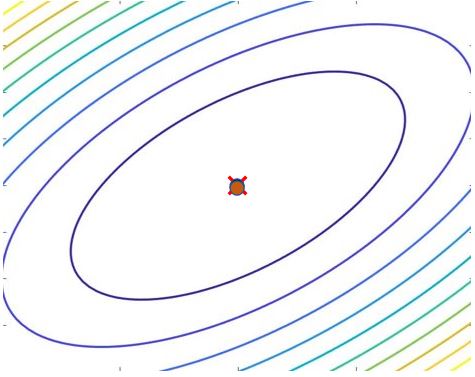}	
	\hspace{0.1cm}
	\includegraphics[scale=0.6]{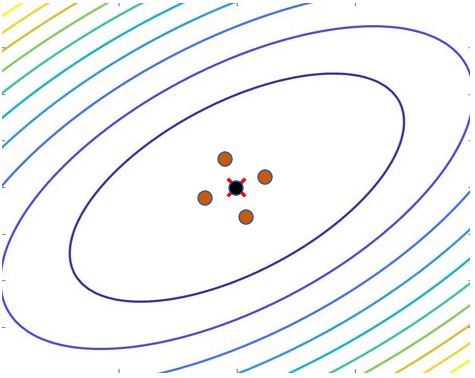}
	\caption{Sketch depicting the minimizers for the case of consensus ($r=0$) (left) and formation control with a single informed agent ($|\mathcal{V}_l|=1$) (right) considered in Lemma \ref{lemm:minimizers}. Black solid circles represents the informed agents, orange solid circles represent the non-informed agents, red cross represents the minimizer of the field.}
	\label{fig:minimizer_sp_cases}	
\end{figure}

Lemma \ref{lemm:minimizers} characterizes the minimizers for the case of consensus ($r=0$) and formation control with a single informed agent ($|\mathcal{V}_l|=1$) as depicted in Fig. \ref{fig:minimizer_sp_cases} and only necessary conditions are obtained for the general case in Lemma \ref{lemm:minimizers_quadratic_radially_symmetric}.
The main tools in deriving these lemmata are definitions of convexity, $\mathcal{L}\mathbf{1}_N=\mathbf{0}$ and the fact that $f$ has a unique minimizer $z$ satisfying $\nabla f(z)=0$.
\begin{lemma}\label{lemm:minimizers}
	Let $\psi$ satisfy Assumption \ref{assum:psi_Sml}, $\mathcal{G}$ and $\mathcal{V}_l$ satisfy Assumption \ref{assm:path_to_leaders} and $f$ be as defined in Definition \ref{defn:f_flock} with $V(y)=\frac{1}{2}(y-r)^T(\mathcal{L} \otimes I_d)(y-r)$.
	Then the following statements hold:
	\begin{enumerate}
		\item[1)] If $r=0$, then $z$ is the minimizer of $f$ if and only if $z=\mathbf{1}_N \otimes y_{\textnormal{opt}}$.
		\item[2)] If $\mathcal{V}_l=\{i\}$ for some $i \in \mathcal{V}$, then $z$ is the minimizer of $f$ if and only if $z_j=y_{\textnormal{opt}}+(r_j-r_i)$ for all $j \in \mathcal{V}$.
	\end{enumerate}
\end{lemma}
\begin{proof}
See Appendix \ref{appendix:proofs}.
\end{proof}
\begin{figure}[b!]	
	\centering
	\includegraphics[scale=0.7]{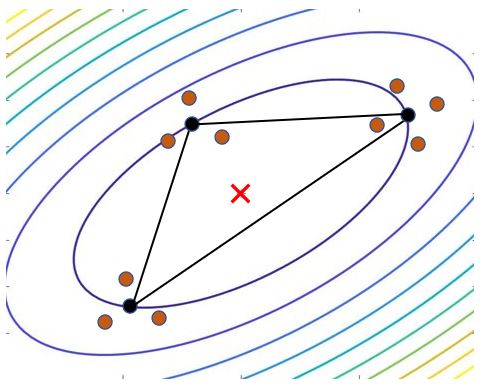}
	\caption{Sketch depicting minimizer for the case of a quadratic field. Black solid circles represents the informed agents, orange solid circles represent the non-informed agents, red cross represents the minimizer of the field as considered in Lemma \ref{lemm:minimizers_quadratic_radially_symmetric} item 2). }
	\label{fig:minimizers_quadratic}	
\end{figure}
\begin{remark}
	Note that if $r$ is chosen such that $r_i=0$ for the informed agent $i \in \mathcal{V}_l$, then $r$ just encodes the desired positions of the agents with the coordinate system such that the source and the informed agent are located at the origin (see Fig. \ref{fig:minimizer_sp_cases} (right)).
\end{remark}
\begin{remark}
	Scenarios involving multiple informed agents and a non-zero formation reference $r$ are difficult to characterize because the terms $\frac{1}{2}(y-r)^T(\mathcal{L} \otimes I_d)(y-r)$ and $\sum_{i \in \mathcal{V}_l} \psi(y_i)$ have competing objectives.
	The equilibrium can thus result in a situation where none of the informed agents are at the source and agents are not at desired relative distances.
\end{remark}

\begin{lemma}\label{lemm:minimizers_quadratic_radially_symmetric}
	Let $\psi$ satisfy Assumption \ref{assum:psi_Sml}, $\mathcal{G}$ and $\mathcal{V}_l$ satisfy Assumption \ref{assm:path_to_leaders}, $f$ be as defined in Definition \ref{defn:f_flock} with $V(y)=\frac{1}{2}(y-r)^T(\mathcal{L} \otimes I_d)(y-r)$ and let $z$ minimize $f$.
	Then $\{y| \nabla \psi(z_i)^T(z_i-y)\geq 0 \textnormal{ for all } i \in \mathcal{V}_l \}$ contains $y_{\textnormal{opt}}$.
	Furthermore, we obtain the following stronger conclusions if $\psi$ is quadratic or radially symmetric.
	\begin{enumerate}
		\item[1)] If $\psi$ is radially symmetric around the source, i.e., it has the form $\psi(y)=\psi_r(||y-y_{\textnormal{opt}}||)$ for some function $\psi_r : \mathbb{R} \rightarrow \mathbb{R}$, then the minimizer $z$ of $f$ is such that $y_{\textnormal{opt}}$ lies in the convex hull of $\{z_i | i \in \mathcal{V}_l\}$. 
		\item[2)] If $\psi$ is quadratic, i.e., it has the form $\psi(y)=y^TQy+b^Ty+c$, then, the minimizer $z$ of $f$ satisfies $y_{\textnormal{opt}}=\frac{1}{|\mathcal{V}_l|} \sum_{i \in \mathcal{V}_l} z_i$, i.e., the center of mass of informed agents is at the minimizer $y_{\textnormal{opt}}$ of $\psi$ (see Fig. \ref{fig:minimizers_quadratic}).
	\end{enumerate}
\end{lemma}
\begin{proof}
See Appendix \ref{appendix:proofs}.
\end{proof}

\begin{figure}[t!]
	\centering 
	\includegraphics[scale=0.4]{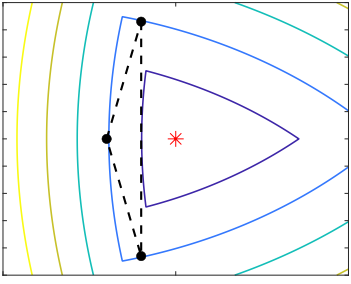}	
	\caption{A sketch of a strongly convex (albeit non-differentiable) scalar field of the form $\psi(y)=\textnormal{max}(c_i^Ty+b_i)+m||y||^2$ with an equilibrium configuration such that the source is outside the convex hull of agent positions (see Remark \ref{rem:counter_example_strongly_convex}).}
	\label{fig:source_seek_minimizers}
\end{figure}
\begin{remark}\label{rem:counter_example_strongly_convex}
	One can construct examples of strongly convex non-quadratic and non-radially symmetric fields such that $y_{\textnormal{opt}}$ does not belong to the convex hull of $\{z_i | i \in \mathcal{V}_l\}$ showing that specializations of quadratic and radially symmetric cases in Lemma \ref{lemm:minimizers_quadratic_radially_symmetric} cannot be extended to the general case without further assumptions on $\psi$.
	See Fig. \ref{rem:counter_example_strongly_convex} for a sketch of a non-smooth strongly convex example field.
\end{remark}

\subsection{Robust Analysis}\label{sec:IQC_theory}
Let us shortly discuss the philosophy behind a typical IQC analysis.
Following the literature on robust control, the uncertain operator $\nabla f$ from dynamics \eqref{eq:sys_dyn_G_hat} is required to belong to a set $\mathbf{\Delta}$ which characterizes our knowledge about the uncertainty in the model.
The first step in carrying out an IQC analysis involves using this knowledge of the set $\mathbf{\Delta}$ to derive properties of the input and output signals $y$ and $u$ in the form of integral inequalities.
These integral inequalities include a factor $e^{2\alpha t}$ (which accounts for the name $\alpha$-IQC) and are specified by a given matrix $P$ and a given LTI system $\Pi=\left[
\begin{array}{c|c}
	A_{\Pi} & B_{\Pi} \\
	\hline
	C_{\Pi} & D_{\Pi}
\end{array}
\right]$ to finally take the form 
\begin{equation}\label{eq:iqc_description}
	\int_0^T e^{2\alpha t}z(t)^T (P \otimes I) z(t) dt \geq 0 \quad \forall T \geq 0,
\end{equation}
where, $z(t)$ is given by
$$
z(t)=\int_0^t C_{\Pi} e^{A_{\Pi}(t-\tau)}B_{\Pi}
\begin{bmatrix}
	y(\tau)\\
	u(\tau)
\end{bmatrix}d\tau
+ 
D_{\Pi}
\begin{bmatrix}
	y(t)\\
	u(t)
\end{bmatrix}.
$$
One aims to derive properties which exactly characterize the set $\mathbf{\Delta}$ by defining an as-large-as-possible class $\mathbb{P}$ such that inequality \eqref{eq:iqc_description} holds for all $P \in \mathbb{P}$.
Once this is done, we can apply standard results such as from \cite{Hu.2016} to get the final analysis result.

Coming back to our problem at hand, we first derive properties between the signals $u(t)$ and $y(t)$ related by the map $u(t)=\nabla f(y(t))$ where $f$ is an arbitrary function in $\mathcal{S}(m,L)$ and $y_*$ minimizes $f$ such that $\nabla f(y_*)=0$.
The standard (non-exponential) ZF IQCs have been well studied.
We give here a purely time-domain derivation of the general non-causal higher order ZF $\alpha$-IQCs.

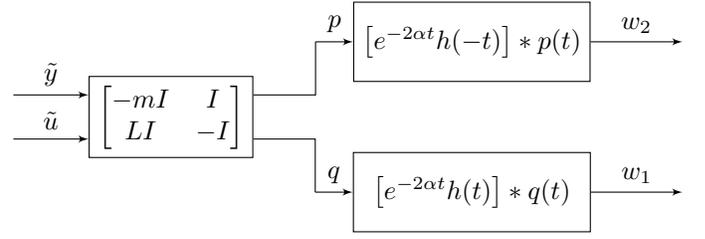
\begin{figure}[t!]
	\centering
	\tikzstyle{block} = [draw, rectangle, 
    minimum height=3em, minimum width=6em]
\tikzstyle{sum} = [draw,  circle, node distance=1cm]
\tikzstyle{input} = [coordinate]
\tikzstyle{output} = [coordinate]
\tikzstyle{pinstyle} = [pin edge={to-,thin,black}]

\begin{tikzpicture}[auto, node distance=2cm,>=latex']
    \node [block,xshift=-1cm] (controller) 
    {$\begin{bmatrix}
    -mI  & I \\
     LI    & -I
    \end{bmatrix}$};
    \node [block, right of=controller,xshift=2cm,yshift=1cm] (filter1) {$\left[e^{-2\alpha t}h(-t)\right]*p(t)$};
    \node [block, right of=controller,xshift=2cm,yshift=-1cm,minimum width=3.15cm] (filter2) {$\left[e^{-2\alpha t}h(t)\right]*q(t)$};
    \node [input, left=of controller.165,xshift=1cm] (input1) {};
    \node [input, left=of controller.195,xshift=1cm] (input2) {};
    \node [input, left=of filter1.west,xshift=1.5cm] (p1) {};
    \node [input, left=of filter2.west,xshift=1.5cm] (q1) {};
    \node [output, right of=filter1,xshift=0.8cm] (output1) {};
    \node [output, right of=filter2,xshift=0.8cm] (output2) {};

    \draw[->] (input1) -- node [name=ytilde] {$\Tilde{y}$}(controller.165);
    \draw[->] (input2) -- node [name=utilde] {$\Tilde{u}$}(controller.195);
    \draw[-] (controller.15) -| node [name=p] {}(p1);
    \draw[-] (controller.345) -| node [name=q] {}(q1);
    \draw[->] (p1) -- node [name=p] {$p$}(filter1.west);
    \draw[->] (q1) -- node [name=q] {$q$}(filter2.west);
    \draw[->] (filter1) -- node [name=w2tilde] {$w_2$}(output1);
    \draw[->] (filter2) -- node [name=w1tilde] {$w_1$}(output2);

\end{tikzpicture}
	\caption{Signal definitions for $p$, $q$, $w_1$, $w_2$.}
	\label{fig:signal_definitions}	
\end{figure}

For convenience, let the deviation signals be defined by
\begin{align*}
	\Tilde{y}(t)&=y(t)-y_*,\\
	\Tilde{u}(t)&=u(t)-u_*=\nabla f(y(t))-\nabla f(y_*)=\nabla f(\Tilde{y}(t)+y_*)
\end{align*}
and for constants $m$ and $L$, let
\begin{equation}\label{eq:signal_defns}
	\begin{split}
		p(t)&=\Tilde{u}(t)-m \Tilde{y}(t),\\
		q(t)&=L \Tilde{y}(t) - \Tilde{u}(t).
	\end{split}
\end{equation}
It can be shown just from basic convexity properties (for example using Proposition 5 from \cite{Lessard.2016}) that if $f \in \mathcal{S}(m,L)$, then $p(t)^Tq(t)\geq 0$ for all $t\geq0$ which implies that $\int_0^Te^{2\alpha t}p(t)^Tq(t) dt \geq 0$ for any $T\geq 0$ and any $\alpha \in \mathbb{R}$.
This corresponds to the well-known sector condition involved in the circle criterion \cite{Scherer.16052021}.
The aim of this section is to derive more general properties which include the sector condition mentioned above in order to better characterize the class of functions $f \in \mathcal{S}(m,L)$.
Let $h \in \mathcal{L}_1(-\infty,\infty)$ be such that for some $H \in \mathbb{R}$
\begin{equation} \label{eq:zf_impulse_resp_cond}
	\begin{aligned}
		h(s) \geq 0 \quad \forall s \in \mathbb{R} \textnormal{   and  }        \int_{-\infty}^{\infty}h(s)ds \leq H,
	\end{aligned}
\end{equation}
 and define signals $w_1(t)$ and $w_2(t)$ by
\begin{equation} \label{eq:signal_defns_w12}
	\begin{split}
		w_1(t)&=\int_0^t e^{-2\alpha (t-\tau)}h(t-\tau) q(\tau) d\tau ,\\
		w_2(t)&=\int_0^t e^{-2\alpha (t-\tau)}h(-(t-\tau)) p(\tau) d\tau .
	\end{split}
\end{equation}
 These signals have been shown in Fig. \ref{fig:signal_definitions} in the form of a block diagram.

\begin{theorem}\label{theom:theorem_pqw_filtered_ZF}
	Let $h$ satisfying \eqref{eq:zf_impulse_resp_cond} be fixed and let $\alpha \geq 0$. 
	Let $\Tilde{u},\Tilde{y} \in \mathcal{L}_{2e}[0,\infty)$ be related by $\Tilde{u}(t)=\nabla f(\Tilde{y}(t)+y_*)$, where $f \in \mathcal{S}(m,L)$ and $y_*$ minimizes $f$.
	Then, the signals defined in \eqref{eq:signal_defns} and \eqref{eq:signal_defns_w12} satisfy the following inequality $\forall T \geq 0$:
	\begin{equation*}\label{eq: IQC_step1}
		\int_0^T e^{2\alpha t}\left(H p(t)^T q(t) - p(t)^T w_1(t) - q(t)^T w_2(t)\right) dt \geq 0.
	\end{equation*}
\end{theorem}
\begin{proof}
	See Appendix \ref{appendix:proofs}.
\end{proof}
\begin{remark}
	A factor of $2H$ is found in the literature instead of the factor $H$ in the first term of the above inequality. 
	This is because $h$ is assumed to be symmetric, i.e., $h(t)=h(-t)$ and is required to satisfy $\int_0^{\infty}h(s)ds\leq H$ (instead of \eqref{eq:zf_impulse_resp_cond}).
\end{remark}
Note that $h(t)\equiv 0$ and $H=1$ satisfies \eqref{eq:zf_impulse_resp_cond} giving us $\int_0^T e^{2\alpha t} p(t)^T q(t) dt \geq 0$, the sector IQC involved in the circle criterion (CC).
Theorem \ref{theom:theorem_pqw_filtered_ZF} thus produces a larger set of IQCs depending on the choice of $h$.
In order to bring the inequality derived in Theorem \ref{theom:theorem_pqw_filtered_ZF} in the typical IQC form \eqref{eq:iqc_description}, we parameterize $h$ satisfying constraints \eqref{eq:zf_impulse_resp_cond}.
This is done along the lines of \cite{Veenman.2016} with suitable modifications and is discussed in Appendix \ref{appendix:h_parameterization}.
The central idea here is to construct $\Pi=\left[
		\begin{array}{c|c}
			A_{\Pi} & B_{\Pi} \\
			\hline
			C_{\Pi} & D_{\Pi}
		\end{array}\right]$ (see Fig. \ref{fig:psi_strucure} from Appendix \ref{appendix:h_parameterization}) and $\mathbb{P}$ such that corresponding to every $P \in \mathbb{P}$, there exists a function $h$ satisfying \eqref{eq:zf_impulse_resp_cond} and the integrand of the inequality in Theorem \ref{theom:theorem_pqw_filtered_ZF} can be written as
	\begin{equation*}
	H p(t)^T q(t) - p(t)^T w_1(t) - q(t)^T w_2(t)=\Tilde{z}^T(t) (P \otimes I) \Tilde{z}(t).
\end{equation*}
where $p(t)$, $q(t)$, $w_1(t)$ and $w_2(t)$ are as defined earlier and 
\begin{equation}\label{eq:signal_defn_z}
	\Tilde{z}(t):=\int_0^t C_{\Pi} e^{A_{\Pi}(t-\tau)}B_{\Pi}
	\begin{bmatrix}
		\Tilde{y}(\tau)\\
		\Tilde{u}(\tau)
	\end{bmatrix}d\tau
	+ 
	D_{\Pi}
	\begin{bmatrix}
		\Tilde{y}(t)\\
		\Tilde{u}(t)
	\end{bmatrix}.
\end{equation}
The next result is thus a direct consequence of Theorem \ref{theom:theorem_pqw_filtered_ZF}.

\begin{theorem}\label{theom:theorem_LMI_ZF}
Let $\Pi$ and $\mathbb{P}$ be as defined in Appendix \ref{appendix:h_parameterization}. and let $\Tilde{u},\Tilde{y} \in \mathcal{L}_{2e}[0,\infty)$ be related by $\Tilde{u}(t)=\nabla f(\Tilde{y}(t)+y_*)$, where $f \in \mathcal{S}(m,L)$ and $y_*$ minimizes $f$.
Then, for any $\alpha \geq 0$, the signal $\Tilde{z}(t)$ as defined in \eqref{eq:signal_defn_z} satisfies
\begin{equation}
    \int_0^T e^{2\alpha t}\Tilde{z}^T(t) (P \otimes I) \Tilde{z}(t) dt \geq 0 \quad \forall P \in \mathbb{P},\forall T \geq 0.
\end{equation}
\end{theorem}
\begin{proof}
See Appendix \ref{appendix:proofs}.
\end{proof}

\begin{remark}\label{rem:CC_ZF_cusal}
	In the classical robust control literature, the transfer function $\hat{\Pi}^*(P \otimes I)\hat{\Pi}$ is the called general non-causal ZF multiplier, where $\hat{\Pi}$ is the transfer function corresponding to the LTI system $\Pi$, $P$ is a matrix belonging to $\mathbb{P}$ and the superscript $*$ represents complex conjugate transpose.
	As discussed in Appendix \ref{appendix:h_parameterization}, the set $\mathbb{P}$ is constructed by variables $H$, $P_1$ and $P_3$.
	When enforcing $P_1=\mathbf{0}$ ($P_3=\mathbf{0}$), the multiplier is called causal (anti-causal) ZF multiplier and when enforcing $P_1=\mathbf{0}$ and $P_3=\mathbf{0}$, one ends up with static multipliers which correspond to the well-known circle-criterion (CC) \cite{Scherer.16052021}.
	The conservatism in these specializations is investigated in Section \ref{sec:numerical_results}.
	Since all the arguments made here are in time-domain, the transfer function $\hat{\Pi}^*(P \otimes I)\hat{\Pi}$ does not appear at any point in this paper.
	The order $\nu$ of the LTI system $\pi$ used for constructing $\Pi$ (see Appendix \ref{appendix:h_parameterization}) is called the order of the ZF multiplier.
\end{remark}
\begin{remark}
	Extension of these results to cases when the map from $\Tilde{y}$ to $\Tilde{u}$ is additionally known to be odd is possible but is not pursued here. 
\end{remark} 
Once we have the IQC result in the form of Theorem \ref{theom:theorem_LMI_ZF} at our disposal, we can apply standard arguments as in \cite{Hu.2016} to obtain the analysis result.
Let $\left[\begin{array}{c|c}
	\mathcal{A}     &  \mathcal{B}\\
	\hline
	\mathcal{C}     &  \mathcal{D}
\end{array}\right]=\Pi \begin{bmatrix}
	I_N \otimes G\\
	I_{Nd}
\end{bmatrix}$ which is used in the next result from \cite{Hu.2016}. 
\begin{theorem}[adapted from \cite{Hu.2016}]\label{theom:main_perf_analysis_formation_big}
	If $\exists \mathcal{X}\succ0,P \in \mathbb{P}$ such that
	\begin{equation}\label{eq:perf_LMI_big}
		\begin{split}
			\begin{bmatrix}
				\mathcal{A}^T\mathcal{X}+\mathcal{X}\mathcal{A}+2\alpha \mathcal{X}  & \mathcal{X}\mathcal{B} \\
				\mathcal{B}^T\mathcal{X}    & \mathbf{0}
			\end{bmatrix}+ 
			\begin{bmatrix}
				\mathcal{C}^T \\ \mathcal{D}^T
			\end{bmatrix}
			P_{(Nd)}
			(*)
			\preceq 
			0,
		\end{split}
	\end{equation}
	then the state trajectories generated by the dynamics \eqref{eq:sys_dyn_G_hat} with any $f \in \mathcal{S}(m,L)$ remain bounded and the output trajectory $y$ converges exponentially to the minimizer $y_*$ of $f$ with rate $\alpha$, i.e., $\exists \kappa \geq 0$ such that  
	$||y(t)-y_*(t)||\leq \kappa e^{-\alpha t}\quad\forall t\geq 0$ .
\end{theorem}
\begin{remark}\label{rem:quasi-LMI}
	Note that although \eqref{eq:perf_LMI_big} is not linear in $\alpha$ and $\mathcal{X}$ due to the product $\alpha \mathcal{X}$, it can be solved efficiently using a bisection over $\alpha$ as commonly suggested in the literature \cite{Lessard.2016}.
\end{remark}
\begin{remark}\label{rem:special_case_single_agent}
	For $N=1$, Theorem \ref{theom:main_perf_analysis_formation_big} reduces to case of single agent in a scalar field considered in \cite{datar2021iqc_source_seek}.
\end{remark}

Extensions of the results obtained in the previous section to LPV systems is straightforward and one such extension is demonstrated next when $N=1$.
Instead of the LTI system $G$, let $G(\rho)$ denote an LPV system with $n_{\rho}$ scheduling parameters \cite{Shamma.1992}, where, for a compact set $\mathcal{P}\subseteq \mathbb{R}^{n_{\rho}} $, the function $\rho:[0,\infty) \rightarrow \mathcal{P}$ captures the time-dependence of the model parameters.
Assume that the system $G(\rho)$ still has the same structure as depicted in Fig. \ref{fig:plant G_flock_chap} with the difference being that the matrices $A,\,B_q,\,B_p,\,C$ may be parameter dependent.
The overall dynamics can be represented by 
\begin{equation} \label{eq:sys_dyn_G_lpv}
	\begin{split}
		\Dot{\eta}(t)&=A_G(\rho(t))\eta(t) + B_G(\rho(t)) u(t), \quad \quad \eta(0)=\eta_0,\\
		y(t)&=C_G(\rho(t)) \eta(t), \\
		u(t)&=\nabla \psi(y(t)),
	\end{split}
\end{equation}
where $\eta\in \mathbb{R}^{n_\eta}$ is the state vector and $\rho:[0,\infty) \rightarrow \mathcal{P}$ is an arbitrary scheduling trajectory.
If the rate of parameter variation $\Dot{\rho}$ is bounded and this bound is known, it could be included by considering parameter-dependent Lyapunov functions (see \cite{scherer_weiland.2000} for details), but this case is not treated here.

Let the series interconnection of the LTI system $(\pi_{m,L} \otimes I_{d})$ and the LPV system $G(\rho)$ be denoted by
$$\left[\begin{array}{c|c}
	\mathcal{A}(\rho)     &  \mathcal{B}(\rho)\\
	\hline
	\mathcal{C}(\rho)     &  \mathcal{D}(\rho)
\end{array}\right] = \Pi \begin{bmatrix}
	G(\rho)\\
	I_{d}
\end{bmatrix}.$$
The following theorem gives a sufficient condition for performance analysis.
\begin{theorem}\label{theom:theorem_perf_analysis_lpv}
If $\exists \mathcal{X}\succ0,P \in \mathbb{P}$ such that, for any $\bar{\rho} \in \mathcal{P}\subseteq \mathbb{R}^{n_{\rho}}$ ($\bar{\rho}$ is a vector, whereas $\rho$ is a function),
\begin{equation}\label{eq:perf_LMI_lpv}
\begin{split}
\begin{bmatrix}
\mathcal{A}(\bar{\rho})^T\mathcal{X}+\mathcal{X}\mathcal{A}(\bar{\rho})+2\alpha \mathcal{X}  & * \\
\mathcal{B}(\bar{\rho})^T\mathcal{X}    & \mathbf{0}
\end{bmatrix}
+
\hspace{0.2cm}
\begin{bmatrix}
	\mathcal{C}(\bar{\rho})^T \\ \mathcal{D}(\bar{\rho})^T 
\end{bmatrix}
P_{(d)}
(*)
\preceq
0,   
\end{split}
\end{equation}
then, the state trajectories generated by the dynamics \eqref{eq:sys_dyn_G_lpv} with any $\psi \in \mathcal{S}(m,L)$ remain bounded and the output trajectory $y$ converges exponentially to the minimizer $y_{\textnormal{opt}}$ of $\psi$ with rate $\alpha$, i.e., $\exists \kappa \geq 0$ such that  
$||y(t)-y_{\textnormal{opt}}||\leq \kappa e^{-\alpha t}$ holds for all $t\geq 0$.
\end{theorem}
\begin{proof}
See Appendix \ref{appendix:proofs}.
\end{proof}

Having demonstrated an extensions to LPV systems, we bring our attention back to the LTI case.
Note that when $N$ is large, LMI \eqref{eq:perf_LMI_big} becomes computationally intractable. 
However, a smaller LMI independent of $N$ can be derived without any additional conservatism.

This is possible due to the specific diagonal and repeated structure of the multiplier and the plant.
The key idea is that once the uncertainty consisting of the interconnections is characterized by an IQC with a diagonal repeated multiplier, the nominal plant and the multiplier form repeated decoupled systems leading to repeated decoupled verification LMIs. See \cite[Section 4.2]{Lessard.2016} for example.

\begin{lemma} \label{lemm:decomposition}
	The following statements are equivalent:
	\begin{enumerate}
		\item $\exists \mathcal{X}\succ 0,P \in \mathbb{P}$ such that \eqref{eq:perf_LMI_big} is satisfied.
		\item $\exists \mathcal{X}_0\succ 0,P \in \mathbb{P}$ such that 
		\begin{equation}\label{eq:perf_LMI}
			\begin{split}
				\begin{bmatrix}
					\mathcal{A}_0^T\mathcal{X}_0+\mathcal{X}_0\mathcal{A}_0+2\alpha \mathcal{X}_0  & * \\
					\mathcal{B}_0^T\mathcal{X}_0    & \mathbf{0}
				\end{bmatrix}
				+
				\begin{bmatrix}
					\mathcal{C}_0^T \\ \mathcal{D}_0^T
				\end{bmatrix}
				P_{(d)}
				(*)
				\preceq
				0,
			\end{split}
		\end{equation}
		where
		$\left[\begin{array}{c|c}
			\mathcal{A}_0     &  \mathcal{B}_0\\
			\hline
			\mathcal{C}_0     &  \mathcal{D}_0
		\end{array}\right]=(\pi_{m,L} \otimes I_d) \begin{bmatrix}
			G\\
			I_{d}
		\end{bmatrix}$
	\end{enumerate}
\end{lemma}
\begin{proof}
See Appendix \ref{appendix:proofs}.
\end{proof}
Finally, let us define an appropriate uncertainty set for dynamics \eqref{eq:sys_dyn_G_hat}.
We assume that $\mathcal{G}$, $\mathcal{V}_l$ and $\psi$ are such that $f \in \mathcal{S}(m,L)$.
This can be verified using Lemma \ref{lemm:f_in_Sml}.
Let $\mathbf{\Delta}_{m,L}$ be a set of tuples $(\mathcal{G},\mathcal{V}_l,\psi)$ such that $f$ as defined in Definition \ref{defn:f_flock} with $V(y)=\frac{1}{2}(y-r)^T(\mathcal{L} \otimes I_d)(y-r)$ belongs to $\mathcal{S}(m,L)$.
\begin{remark}
	Note that with $\mathcal{V}=\mathcal{V}_l=\{1\}$ and $\mathcal{E}$ empty, we obtain $f=\psi$ and the setup reduces to the scenario of a single agent embedded in a scalar field considered in \cite{datar2021iqc_source_seek}. 
\end{remark}
\begin{theorem} \label{theom:main_perf_analysis_formation}
	Let graph $\mathcal{G}$, set of informed agents $\mathcal{V}_l$ and scalar field $\psi$ be such that $(\mathcal{G},\mathcal{V}_l,\psi) \in \mathbf{\Delta}_{m,L}$ for some $0<m\leq L$.
	Let $y_*$ be the minimizer of $f$.
	If $\exists \mathcal{X}_0>0,P \in \mathbb{P}$ such that \eqref{eq:perf_LMI} is satisfied, then, under dynamics \eqref{eq:sys_dyn_G_hat}, $y$ converges exponentially to $y_*$ with rate $\alpha$, i.e., $\exists \kappa \geq 0$ such that,   
	$||y(t)-y_*||\leq \kappa e^{-\alpha t}$ holds for all $t\geq 0$.
\end{theorem}
\begin{proof}
	Since $(\mathcal{G},\mathcal{V}_l,\psi) \in \mathbf{\Delta}_{m,L}$, $f \in \mathcal{S}(m,L)$.
	Using the hypothesis of this theorem and Lemma \ref{lemm:decomposition}, $\exists \mathcal{X}\succ 0,P \in \mathbb{P}$ such that \eqref{eq:perf_LMI_big} is satisfied.
	Finally, application of Theorem \ref{theom:main_perf_analysis_formation_big} completes the proof. 
\end{proof}
\begin{remark}
	These results can be easily extended to the setting of Linear Parameter Varying (LPV) systems or even non-linear systems with a quasi-LPV representation.
	Note that although this would lead to heterogeneous scheduling in the networked LPV system, the decomposition result presented in Lemma \ref{lemm:decomposition} can still be extended to such quasi-LPV systems as long as the parameter sets $\mathcal{P}$ is the same for all agents (see \cite{datar2023gradient}).
\end{remark}

\section{Analysis: Flocking (Problem 2)}
The analysis for Problem \ref{prob:flocking} is based on standard dissipativity based arguments.
To exploit the IQC descriptions from Assumption \ref{assum:f_flock}, items $3)$ and $4)$, dynamics \eqref{eq:sys_dyn_G_hat} are transformed into a suitable form to facilitate the analysis.
Specifically, we write the input as
$\nabla f(y(t))=\nabla f(q(t))-(\nabla f(q(t))-\nabla f(y(t)))$
and consider as output channels $q(t)$ and  $q(t)-y(t)$.
Additionally, the state $\eta=[\eta_1^T\;\cdots\;\eta_N^T]$ is permuted to a new state vector $\bar{\eta}=[x^T\;q^T\;p^T]^T$ to obtain equivalent dynamics

\begin{equation} \label{eq:sys_dyn_flock_equival}
	\begin{split}
		\begin{bmatrix}
			\Dot{\bar{\eta}}(t)\\
			{q}(t)\\
			{q}(t)-{y}(t)
		\end{bmatrix}
		&=
		\left[\begin{array}{ccc}
			\mathcal{A}_G     &  \mathcal{B}_G&  -\mathcal{B}_G\\
			\mathcal{C}_{G1}     &  \mathbf{0}  &  \mathbf{0}\\
			\mathcal{C}_{G2}     &  \mathbf{0} &  \mathbf{0}\\
		\end{array}\right]
		\begin{bmatrix}
			\bar{\eta}(t)\\
			{d}_1(t)\\
			{d}_2(t)
		\end{bmatrix},\\
		{d}_1(t)&=\nabla f ({q}(t)), \\
		{d}_2(t)&=\nabla f ({q}(t))-\nabla f ({y}(t)),
	\end{split}
\end{equation}
where the model matrices are given in Appendix \ref{appendix:model_matrices}.

Drawing motivation from \cite[Section 3]{Fazlyab.2018}, a non-negative function $V_s$ of the state $\bar{\eta}$ is constructed such that it is $0$ only when $q_*=y_*$, $x_*=-A^{-1}B_q y_*$ and $p_*=0$ where $y_* \in \mathcal{S}$ is the local minimizer as per Assumption \ref{assum:f_flock}.
In this regard, for given matrices $A$, $B_q$ and matrix variables $R\succeq 0$ and $Q\succ0$, construct a block $3 \times 3$ matrix  $\mathcal{X}_0$ as 
\begin{align*}\label{eq:cal_X_0}
	\begin{bmatrix}
		\mathcal{X}_{11}    & \mathcal{X}_{12} & \mathcal{X}_{13} \\
		\mathcal{X}_{21}    & \mathcal{X}_{22} & \mathcal{X}_{23} \\
		\mathcal{X}_{31}    & \mathcal{X}_{32} & \mathcal{X}_{33} \\
	\end{bmatrix}&= R +
	\begin{bmatrix}
		Q    & QA^{-1}B_q & \mathbf{0} \\
		*    & B_q^TA^{-T}QA^{-1}B_q   & \mathbf{0} \\
		\mathbf{0}    & \mathbf{0} & I_d
	\end{bmatrix}. 
\end{align*}
With $\bar{\eta}_*=[x_*^T\;q_*^T\;p_*^T]^T$, define the storage function 
\begin{equation*} \label{eq: V_flock}
	V_s(x,q,p)= 
	(*)
	\begin{bmatrix}
		\hat{\mathcal{X}}_{11}    & \hat{\mathcal{X}}_{12} & \hat{\mathcal{X}}_{13} \\
		\hat{\mathcal{X}}_{21}    & \hat{\mathcal{X}}_{22} & \hat{\mathcal{X}}_{23} \\
		\hat{\mathcal{X}}_{31}    & \hat{\mathcal{X}}_{32} & \hat{\mathcal{X}}_{33} \\
	\end{bmatrix}
	(\bar{\eta}-\bar{\eta}_*)
	+2\mu (f(q)-f_{\textnormal{min}}).
\end{equation*}

With an energy function $E(t):=V_s(x(t),q(t),p(t))$, it can be shown that if $E(t)$ is below some prescribed value, the state is bounded and the output belongs to $\mathcal{S}_0\subset\mathcal{S}$.
Feasibility of $\mathcal{Z}\preceq 0$ from \eqref{eqn:LMI_flock_Z} implies that if the output is in $\mathcal{S}$, $E(t)$ is non-increasing.
If the initial conditions are such that $E(0)$ is small enough and initial output is in $\mathcal{S}_0$ (ensured by \eqref{eq:ic_small}), a forward invariance of set $\mathcal{S}_0$ and thus boundedness of state and output trajectories is established.
LaSalle's invariance principle then implies convergence.
This is summarized in the next theorem.
\begin{figure*}[t!]
	\normalsize
	
	
	\begin{equation}\label{eqn:LMI_flock_Z}
		\mathcal{Z}=\begin{bmatrix}
			\mathcal{A}_0^T\mathcal{X}_0+\mathcal{X}_0\mathcal{A}_0  & \mathcal{X}_0\mathcal{B}_0 \\
			\mathcal{B}_0^T\mathcal{X}_0    & \mathbf{0}
		\end{bmatrix}
		+
		\begin{bmatrix}
			\mathbf{0} & \mathbf{0}& \mathbf{0}& \mathbf{0}& \mathbf{0} \\
			\mathbf{0} & \mathbf{0}& \mathbf{0}& \mathbf{0}& \mathbf{0} \\
			\mathbf{0} & \mathbf{0}& \varepsilon I_d& \mu I_d&      \mathbf{0} \\
			\mathbf{0} & \mathbf{0}& \mu I_d&      \mathbf{0}& \mathbf{0} \\
			\mathbf{0} & \mathbf{0}& \mathbf{0}& \mathbf{0}& \mathbf{0} \\
		\end{bmatrix}
		+
		\sum_{i=1}^2
		(*)
		(\lambda_i M_{i0} \otimes I_d)
		\begin{bmatrix}
			\mathcal{C}_{i0} & \mathcal{D}_{i0}\\
		\end{bmatrix}
	\end{equation}
	
	\hrulefill
	\vspace*{2pt}

\end{figure*}
\begin{theorem}[Analysis for Problem \ref{prob:flocking}] \label{theom:flocking_analysis_local}
    Let $f: \mathbb{R}^{Nd} \rightarrow \mathbb{R}$ and an open set $\mathcal{S} \subseteq \mathbb{R}^{Nd}$ satisfy Assumption \ref{assum:f_flock} and let the vehicle dynamics \eqref{eq:vehicle_dyn_flock_chap} satisfy Assumption \ref{assm:ss_error}.
    If there exist $R\succeq 0$, $Q\succ 0$, $\mu>0$, $\lambda_1\geq0$, $\lambda_2\geq0$ and $\varepsilon>0$ such that $\mathcal{Z}$ defined in \eqref{eqn:LMI_flock_Z} with model matrices given in Appendix \ref{appendix:model_matrices} is negative semi-definite
    and if the initial conditions $x_0$, $q_0$, $p_0$ are such that 
    \begin{equation}\label{eq:ic_small}
    	V_s(x_0,q_0,p_0)\leq \textnormal{min}\{2 c_1 \mu , \frac{{ c_2\lambda_{\textnormal{min}}(Q)}}{||C||^2} \},
    \end{equation}
    where $c_1$ and $c_2$ are constants involved in Assumption \ref{assum:f_flock} item 5),
    then the state trajectory $\eta$ generated by the dynamics \eqref{eq:sys_dyn_G_hat} remains bounded for all $t\geq 0$ and $y(t)$ converges to the set $\{y_*\in \mathcal{S}|\nabla f(y_*)=0\}$.
\end{theorem}
\begin{proof}
See Appendix \ref{appendix:proofs}.
\end{proof}
\begin{remark}
	Note that the sufficient condition $\mathcal{Z}\preceq 0$ (see \eqref{eqn:LMI_flock_Z}) contains model matrices appearing in the dynamics of a single agent and is independent of the network-size $N$.
\end{remark}

\begin{corollary}\label{corr:global_flocking_analysis}
	Let $f$, $\mathcal{S}=\mathbb{R}^{Nd}$ satisfy Assumption \ref{assum:f_flock} and assume that the vehicle dynamics \eqref{eq:vehicle_dyn_flock_chap} satisfy Assumption \ref{assm:ss_error}.
	Additionally, let $f$ be proper, i.e., $\{q|f(q)\leq c\}$ is compact for all $c \in \mathbb{R}$ \cite{FB-LNS}. 
	If there exist $R\succeq 0$, $Q\succ 0$, $\mu>0$, $\lambda_1\geq0$, $\lambda_2\geq0$ such that $\mathcal{Z}$ defined in \eqref{eqn:LMI_flock_Z} with model matrices given in Appendix \ref{appendix:model_matrices} is negative semi-definite, then, for any initial condition, the state trajectory $\eta$ under dynamics \eqref{eq:sys_dyn_G_hat} remains bounded for all $t\geq 0$ and $y(t)$ converges to the set $\{ y_*|\nabla f(y_*)=0\}$.
\end{corollary}
\begin{proof}
	Since $f$ is proper and $\mathcal{S}=\mathbb{R}^{Nd}$, the set $\mathcal{S}_0$ is bounded and contained in $\mathcal{S}$ for any $c_1$, $c_2$.
	Therefore, the requirement on the initial condition is satisfied.
	Applying Theorem \ref{theom:flocking_analysis_local} with the above observation gives the result.
\end{proof}
\begin{remark}
	Theorem \ref{theom:flocking_analysis_local} and Corollary \ref{corr:global_flocking_analysis} are presented under general assumptions on $f$ without specifying the $M_{10}$ and $M_{20}$.
	Examples use $M_{10}=M_{20}=\begin{bmatrix}
		L^2&0\\0&-1
	\end{bmatrix}$. 
\end{remark}

The code used for generating results in the following Sections \ref{sec:numerical_results}, \ref{sec:numerical_results_formation_control} and \ref{sec:numerical_results_flocking} is availabe at \cite{datar_2022_6672230}.
\section{Numerical Examples: Single Agent}\label{sec:numerical_results}

\subsection{Single quadrotor in $\psi \in \mathcal{S}(1,L_{\psi})$}\label{eg:quadrotor}

We consider a linearized quadrotor model and use an LQR (Linear-Quadratic-Regulator) based state-feedback controller tuned for zero steady-state error for step references.  
As discussed in Section \ref{sec:problem_setup}, we let this closed-loop system be represented by state-space realization \eqref{eq:vehicle_dyn_flock_chap} and augment it with dynamics \eqref{eq:vir_vehicle_dyn_flock_chap} to obtain the system $G$. 

\begin{figure}[t!]
	\centering
%
%
\definecolor{mycolor1}{rgb}{0.00000,1.00000,1.00000}%
\begin{tikzpicture}

\begin{axis}[%
width=2.521in,
height=2.366in,
at={(0.758in,0.481in)},
scale only axis,
xmin=1,
xmax=9,
xlabel style={font=\color{white!15!black}},
xlabel={$L_{\psi}$},
ymin=0,
ymax=0.5,
ylabel style={font=\color{white!15!black}},
ylabel={$\alpha$},
axis background/.style={fill=white},
title style={font=\bfseries},
title={},
legend style={legend cell align=left, align=left, draw=white!15!black}
]
\addplot [color=red, line width=1.0pt, draw=none, mark=o, mark options={solid, red}]
  table[row sep=crcr]{%
1	0.4119873046875\\
1.2	0.3472900390625\\
1.4	0.3076171875\\
1.6	0.277099609375\\
1.8	0.2508544921875\\
2	0.2276611328125\\
2.2	0.2069091796875\\
2.4	0.1873779296875\\
2.6	0.1690673828125\\
2.8	0.152587890625\\
3	0.1361083984375\\
3.2	0.120849609375\\
3.4	0.106201171875\\
3.6	0.0921630859375\\
3.8	0.0787353515625\\
4	0.06591796875\\
4.2	0.0531005859375\\
4.4	0.0408935546875\\
4.6	0.029296875\\
4.8	0.0177001953125\\
5	0.006103515625\\
5.2	-1\\
5.4	-1\\
5.6	-1\\
5.8	-1\\
6	-1\\
6.2	-1\\
6.4	-1\\
6.6	-1\\
6.8	-1\\
7	-1\\
7.2	-1\\
7.4	-1\\
7.6	-1\\
7.8	-1\\
8	-1\\
8.2	-1\\
8.4	-1\\
8.6	-1\\
8.8	-1\\
9	-1\\
};
\addlegendentry{CC (static)}

\addplot [color=green, dashed, line width=1.0pt]
  table[row sep=crcr]{%
1	0.4119873046875\\
1.2	0.3643798828125\\
1.4	0.33935546875\\
1.6	0.321044921875\\
1.8	0.306396484375\\
2	0.2935791015625\\
2.2	0.281982421875\\
2.4	0.27099609375\\
2.6	0.260009765625\\
2.8	0.2490234375\\
3	0.2386474609375\\
3.2	0.2276611328125\\
3.4	0.21728515625\\
3.6	0.2069091796875\\
3.8	0.196533203125\\
4	0.1861572265625\\
4.2	0.17578125\\
4.4	0.1654052734375\\
4.6	0.1556396484375\\
4.8	0.145263671875\\
5	0.135498046875\\
5.2	0.125732421875\\
5.4	0.115966796875\\
5.6	0.1068115234375\\
5.8	0.0970458984375\\
6	0.087890625\\
6.2	0.0787353515625\\
6.4	0.069580078125\\
6.6	0.06103515625\\
6.8	0.0518798828125\\
7	0.0433349609375\\
7.2	0.0347900390625\\
7.4	0.0262451171875\\
7.6	0.0177001953125\\
7.8	0.009765625\\
8	0.0018310546875\\
8.2	-1\\
8.4	-1\\
8.6	-1\\
8.8	-1\\
9	-1\\
};
\addlegendentry{ZF causal order 1}

\addplot [color=blue, dashdotted, line width=1.0pt]
  table[row sep=crcr]{%
1	0.4119873046875\\
1.2	0.3472900390625\\
1.4	0.3076171875\\
1.6	0.277099609375\\
1.8	0.2508544921875\\
2	0.2276611328125\\
2.2	0.2069091796875\\
2.4	0.1873779296875\\
2.6	0.1690673828125\\
2.8	0.152587890625\\
3	0.1361083984375\\
3.2	0.120849609375\\
3.4	0.106201171875\\
3.6	0.0921630859375\\
3.8	0.0787353515625\\
4	0.06591796875\\
4.2	0.0531005859375\\
4.4	0.0408935546875\\
4.6	0.029296875\\
4.8	0.0177001953125\\
5	0.006103515625\\
5.2	-1\\
5.4	-1\\
5.6	-1\\
5.8	-1\\
6	-1\\
6.2	-1\\
6.4	-1\\
6.6	-1\\
6.8	-1\\
7	-1\\
7.2	-1\\
7.4	-1\\
7.6	-1\\
7.8	-1\\
8	-1\\
8.2	-1\\
8.4	-1\\
8.6	-1\\
8.8	-1\\
9	-1\\
};
\addlegendentry{ZF anti-causal order 1}

\addplot [color=mycolor1, line width=1.0pt]
  table[row sep=crcr]{%
1	0.4119873046875\\
1.2	0.3955078125\\
1.4	0.379638671875\\
1.6	0.3643798828125\\
1.8	0.34912109375\\
2	0.3338623046875\\
2.2	0.3192138671875\\
2.4	0.30517578125\\
2.6	0.2911376953125\\
2.8	0.2777099609375\\
3	0.2642822265625\\
3.2	0.25146484375\\
3.4	0.2386474609375\\
3.6	0.2264404296875\\
3.8	0.2142333984375\\
4	0.20263671875\\
4.2	0.1910400390625\\
4.4	0.179443359375\\
4.6	0.16845703125\\
4.8	0.157470703125\\
5	0.146484375\\
5.2	0.1361083984375\\
5.4	0.125732421875\\
5.6	0.1153564453125\\
5.8	0.1055908203125\\
6	0.0958251953125\\
6.2	0.0860595703125\\
6.4	0.0762939453125\\
6.6	0.067138671875\\
6.8	0.0579833984375\\
7	0.048828125\\
7.2	0.040283203125\\
7.4	0.0311279296875\\
7.6	0.0225830078125\\
7.8	0.0140380859375\\
8	0.006103515625\\
8.2	-1\\
8.4	-1\\
8.6	-1\\
8.8	-1\\
9	-1\\
};
\addlegendentry{ZF order 1}

\addplot [color=black, line width=1.0pt, draw=none, mark=asterisk, mark options={solid, black}]
  table[row sep=crcr]{%
1	0.412331651029501\\
1.2	0.395974289263064\\
1.4	0.380000019452199\\
1.6	0.364400808966183\\
1.8	0.349167588131177\\
2	0.334290538830461\\
2.2	0.319759335099553\\
2.4	0.305563341010742\\
2.6	0.291691771302537\\
2.8	0.27813382002905\\
3	0.264878762120048\\
3.2	0.251916032253449\\
3.4	0.239235284914788\\
3.6	0.226826438996403\\
3.8	0.2146797097987\\
4	0.202785630850186\\
4.2	0.191135067568535\\
4.4	0.179719224441611\\
4.6	0.168529647113125\\
4.8	0.157558220508057\\
5	0.146797163923227\\
5.2	0.136239023833451\\
5.4	0.125876665018557\\
5.6	0.115703260496915\\
5.8	0.105712280652737\\
6	0.095897481864014\\
6.2	0.08625289487244\\
6.4	0.0767728130832965\\
6.6	0.067451780940349\\
6.8	0.0582845824858749\\
7	0.0492662301880655\\
7.2	0.0403919540956929\\
7.4	0.0316571913620181\\
7.6	0.0230575761658658\\
7.8	0.0145889300465953\\
8	0.00624725266099713\\
8.2	-0.00197128703663169\\
8.4	-0.0100703591950899\\
8.6	-0.0180534810554159\\
8.8	-0.0259240245392554\\
9	-0.0336852234492348\\
};
\addlegendentry{Example fields $\psi_0$}

\end{axis}
\end{tikzpicture}%
	\caption{Convergence rate estimates for quadrotor dynamics provided by different multipliers (see Remark \ref{rem:CC_ZF_cusal}) for fields $\psi\in \mathcal{S}(1,L_{\psi})$ and convergence rates for $\psi_0(x)=\frac{1}{2}(x-y_{\textnormal{opt}})^T\textnormal{diag}(1,L_{\psi})(x-y_{\textnormal{opt}})$}
	\label{fig:quadrotor_robustness}
\end{figure}
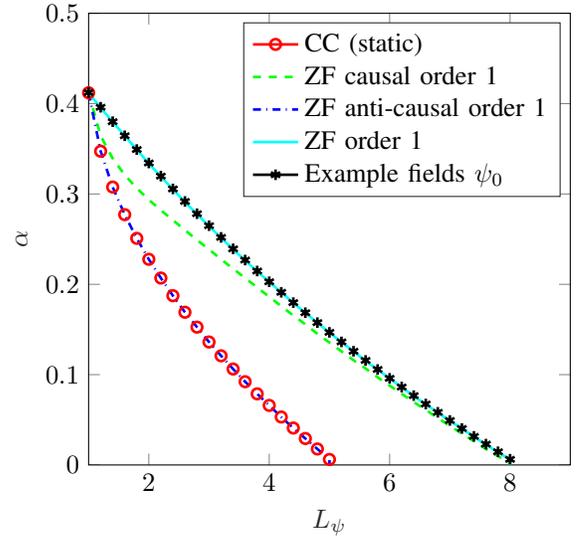

For fixed gains $k_p$ and $k_d$ and given closed-loop quadrotor dynamics, Fig. \ref{fig:quadrotor_robustness} shows the convergence rate estimates provided by different multipliers for fields $\psi\in \mathcal{S}(1,L_{\psi})$ with increasing $L_{\psi}$.
Since increasing $L_{\psi}$ enlarges the set of allowable fields, i.e., $\mathcal{S}(1,L_1) \subseteq \mathcal{S}(1,L_2) \quad \forall L_1\leq L_2$, the estimates are non-increasing with increasing $L_{\psi}$. 
It can be seen that while we can certify stability with the circle criterion for fields $\psi \in \mathcal{S}(1,5.1)$, the general non-causal ZF multipliers along with the ZF multipliers restricted to the causal case ($P_3=0$) can certify stability for all fields $\psi \in \mathcal{L}(1,8.1)$.
Furthermore, for each $L_{\psi}$, $\psi(x)=\frac{1}{2}(x-y_{\textnormal{opt}})^T\begin{bmatrix}1&\\&L_{\psi}\end{bmatrix}(x-y_{\textnormal{opt}})$ achieves the convergence rate guaranteed by the analysis showing that, in this example, the estimates are tight.
The conservatism incurred by restricting the search to causal multipliers is minor in this example.
The stability analysis in \cite{Attallah.2020} uses manually constructed diagonal storage functions together with a small gain argument and for this example, gives the sufficient condition for stability to be $L_{\psi}<5$.
This interestingly coincides with the stability boundary given by the circle-criterion (static multipliers).
Since performance analysis was not included in \cite{Attallah.2020}, this example illustrates the extension of \cite{Attallah.2020} to a non-conservative performance analysis.

\subsection{Example showing the benefit of non-causal multipliers}
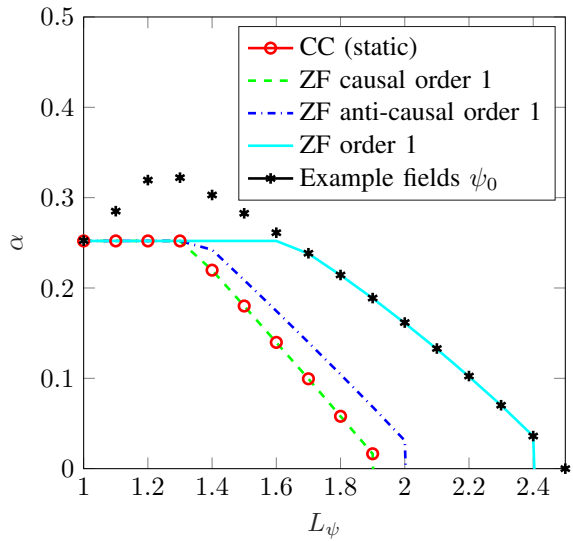
\begin{figure}[t!]
	\centering
%
%
\definecolor{mycolor1}{rgb}{0.00000,1.00000,1.00000}%
\begin{tikzpicture}

\begin{axis}[%
width=2.521in,
height=2.366in,
at={(0.758in,0.481in)},
scale only axis,
xmin=1,
xmax=2.5,
xlabel style={font=\color{white!15!black}},
xlabel={$L_{\psi}$},
ymin=0,
ymax=0.5,
ylabel style={font=\color{white!15!black}},
ylabel={$\alpha$},
axis background/.style={fill=white},
title style={font=\bfseries},
title={},
legend style={legend cell align=left, align=left, draw=white!15!black}
]
\addplot [color=red, line width=1.0pt, draw=none, mark=o, mark options={solid, red}]
  table[row sep=crcr]{%
1	0.2520751953125\\
1.1	0.2520751953125\\
1.2	0.2520751953125\\
1.3	0.2520751953125\\
1.4	0.2197265625\\
1.5	0.1800537109375\\
1.6	0.1397705078125\\
1.7	0.0994873046875\\
1.8	0.0579833984375\\
1.9	0.0164794921875\\
2	-1\\
2.1	-1\\
2.2	-1\\
2.3	-1\\
2.4	-1\\
2.5	-1\\
2.6	-1\\
2.7	-1\\
2.8	-1\\
2.9	-1\\
3	-1\\
};
\addlegendentry{CC (static)}

\addplot [color=green, dashed, line width=1.0pt]
  table[row sep=crcr]{%
1	0.2520751953125\\
1.1	0.2520751953125\\
1.2	0.2520751953125\\
1.3	0.2520751953125\\
1.4	0.2197265625\\
1.5	0.1800537109375\\
1.6	0.1397705078125\\
1.7	0.0994873046875\\
1.8	0.0579833984375\\
1.9	0.0164794921875\\
2	-1\\
2.1	-1\\
2.2	-1\\
2.3	-1\\
2.4	-1\\
2.5	-1\\
2.6	-1\\
2.7	-1\\
2.8	-1\\
2.9	-1\\
3	-1\\
};
\addlegendentry{ZF causal order 1}

\addplot [color=blue, dashdotted, line width=1.0pt]
  table[row sep=crcr]{%
1	0.2520751953125\\
1.1	0.2520751953125\\
1.2	0.2520751953125\\
1.3	0.2520751953125\\
1.4	0.2423095703125\\
1.5	0.208740234375\\
1.6	0.174560546875\\
1.7	0.1397705078125\\
1.8	0.1043701171875\\
1.9	0.068359375\\
2	0.0311279296875\\
2.1	-1\\
2.2	-1\\
2.3	-1\\
2.4	-1\\
2.5	-1\\
2.6	-1\\
2.7	-1\\
2.8	-1\\
2.9	-1\\
3	-1\\
};
\addlegendentry{ZF anti-causal order 1}

\addplot [color=mycolor1, line width=1.0pt]
  table[row sep=crcr]{%
1	0.2520751953125\\
1.1	0.2520751953125\\
1.2	0.2520751953125\\
1.3	0.2520751953125\\
1.4	0.2520751953125\\
1.5	0.2520751953125\\
1.6	0.2520751953125\\
1.7	0.238037109375\\
1.8	0.2142333984375\\
1.9	0.1885986328125\\
2	0.1611328125\\
2.1	0.1324462890625\\
2.2	0.1019287109375\\
2.3	0.0701904296875\\
2.4	0.0360107421875\\
2.5	-1\\
2.6	-1\\
2.7	-1\\
2.8	-1\\
2.9	-1\\
3	-1\\
};
\addlegendentry{ZF order 1}

\addplot [color=black, line width=1.0pt, draw=none, mark=asterisk, mark options={solid, black}]
  table[row sep=crcr]{%
1	0.252381023686084\\
1.1	0.285030034676752\\
1.2	0.319444587967699\\
1.3	0.322120595821607\\
1.4	0.30294136645159\\
1.5	0.282662461122388\\
1.6	0.261201184778797\\
1.7	0.238470096600161\\
1.8	0.214377762177947\\
1.9	0.188829987346274\\
2	0.161731674331652\\
2.1	0.132989449419385\\
2.2	0.102515199007407\\
2.3	0.0702306031940527\\
2.4	0.0360726564236493\\
2.5	-1.11022302462516e-16\\
2.6	-0.0380003387049316\\
2.7	-0.0779064749243448\\
2.8	-0.119656383851905\\
2.9	-0.163144853458106\\
3	-0.208223620537594\\
};
\addlegendentry{Example fields $\psi_0$}

\end{axis}
\end{tikzpicture}%
	\caption{Convergence rate estimates for $G(s)=\frac{5(s-1)}{s(s^2+s+25)}$ provided by different multipliers (see Remark \ref{rem:CC_ZF_cusal}) for fields $\psi\in \mathcal{S}(1,L_{\psi})$ and convergence rates for $\psi_0(x)=\frac{1}{2} L_{\psi}x^2$} 
	\label{fig:robustness_non_causal_multiplier}
\end{figure}

We now present an academic example that brings out the benefit of using general non-causal multipliers over causal multipliers.
Let $G(s)=5\frac{(s-1)}{s(s^2+s+25)}$ and consider fields $\psi\in \mathcal{S}(1,L_{\psi})$.
The convergence rate estimates provided by different multipliers for increasing $L$ is shown in Fig. \ref{fig:robustness_non_causal_multiplier}.
It can be seen that while the circle criterion and causal ZF multipliers certify stability for fields $\psi\in \mathcal{S}(1,1.9)$, the anti-causal ZF multipliers can certify stability for fields $\psi\in \mathcal{S}(1,2)$ and the general non-causal ZF multipliers can certify stability for fields $\psi\in \mathcal{S}(1,2.4)$.
Furthermore, convergence rates with $\psi_0(x)=\frac{1}{2} L_{\psi}x^2$ hit the convergence rate estimates showing that these estimates are tight.
Note that the gap between the actual convergence rates for example fields and the estimates obtained from non-causal ZF multipliers for $L_{\psi}\in \{ 1.1, 1.2, \cdots, 1.6 \} $ does not imply conservatism since the example field used at $L_{\psi}=1$ is included for any larger $L_{\psi}$.

\subsection{LPV generic vehicle model}
\begin{figure}[t!]
    	\centering
%
%
\definecolor{mycolor1}{rgb}{0.00000,1.00000,1.00000}%
\definecolor{mycolor2}{rgb}{0.00000,0.44700,0.74100}%
\definecolor{mycolor3}{rgb}{0.85000,0.32500,0.09800}%
\definecolor{mycolor4}{rgb}{0.92900,0.69400,0.12500}%
\definecolor{mycolor5}{rgb}{0.49400,0.18400,0.55600}%
\begin{tikzpicture}

\begin{axis}[%
width=2.521in,
height=2.366in,
at={(0.758in,0.481in)},
scale only axis,
xmin=0,
xmax=30,
xlabel style={font=\color{white!15!black}},
xlabel={$L_{\psi}$},
ymin=0,
ymax=0.5,
ylabel style={font=\color{white!15!black}},
ylabel={$\alpha$},
axis background/.style={fill=white},
title style={font=\bfseries},
title={},
axis x line*=bottom,
axis y line*=left,
legend style={legend cell align=left, align=left, draw=white!15!black},
legend pos=south east
]
\addplot [color=green, dashed, line width=1.0pt]
  table[row sep=crcr]{%
1	0.3900146484375\\
3	0.1971435546875\\
5	0.0567626953125\\
7	-1\\
9	-1\\
11	-1\\
13	-1\\
15	-1\\
17	-1\\
19	-1\\
21	-1\\
23	-1\\
25	-1\\
27	-1\\
29	-1\\
};
\addlegendentry{ZF casual order 1}


\addplot [color=mycolor1, line width=1.0pt]
  table[row sep=crcr]{%
1	0.3900146484375\\
3	0.36376953125\\
5	0.3253173828125\\
7	0.29052734375\\
9	0.2618408203125\\
11	0.2374267578125\\
13	0.2166748046875\\
15	0.198974609375\\
17	0.1837158203125\\
19	0.1702880859375\\
21	0.15869140625\\
23	0.1483154296875\\
25	0.13916015625\\
27	0.130615234375\\
29	0.123291015625\\
};
\addlegendentry{ZF order 1}

\addplot [color=mycolor2, dotted, line width=1.0pt]
  table[row sep=crcr]{%
1	0.3900146484375\\
3	0.3900146484375\\
5	0.3900146484375\\
7	0.38818359375\\
9	0.3765869140625\\
11	0.3607177734375\\
13	0.3448486328125\\
15	0.3289794921875\\
17	0.3143310546875\\
19	0.2996826171875\\
21	0.2862548828125\\
23	0.2740478515625\\
25	0.2618408203125\\
27	0.2508544921875\\
29	0.2398681640625\\
};
\addlegendentry{ZF order 2}



\addplot [color=mycolor5, dotted, line width=1.0pt]
  table[row sep=crcr]{%
1	0.3900146484375\\
3	0.3900146484375\\
5	0.3900146484375\\
7	0.3900146484375\\
9	0.3900146484375\\
11	0.3900146484375\\
13	0.3900146484375\\
15	0.3900146484375\\
17	0.3900146484375\\
19	0.3900146484375\\
21	0.3900146484375\\
23	0.3900146484375\\
25	0.3900146484375\\
27	0.3900146484375\\
29	0.3900146484375\\
};
\addlegendentry{ZF order 5}

\addplot [color=black, line width=1.0pt, draw=none, mark=asterisk, mark options={solid, black}]
  table[row sep=crcr]{%
1	0.4\\
3	0.4\\
5	0.4\\
7	0.4\\
9	0.4\\
11	0.4\\
13	0.4\\
15	0.4\\
17	0.4\\
19	0.4\\
21	0.4\\
23	0.4\\
25	0.4\\
27	0.4\\
29	0.4\\
};
\addlegendentry{Example fields $\psi_0$}

\end{axis}
\end{tikzpicture}%
    	\caption{Convergence rate estimates for LPV system \eqref{eq:lpv_example} provided by different multipliers (see Remark \ref{rem:CC_ZF_cusal}) for fields $\psi\in \mathcal{S}(1,L_{\psi})$ and convergence rates for $\psi_0(x)=\frac{1}{2}(x-y_{\textnormal{opt}})^T\textnormal{diag}(1,L_{\psi})(x-y_{\textnormal{opt}})$} 
    	\label{fig:LPV_generic_vehicle}
\end{figure}
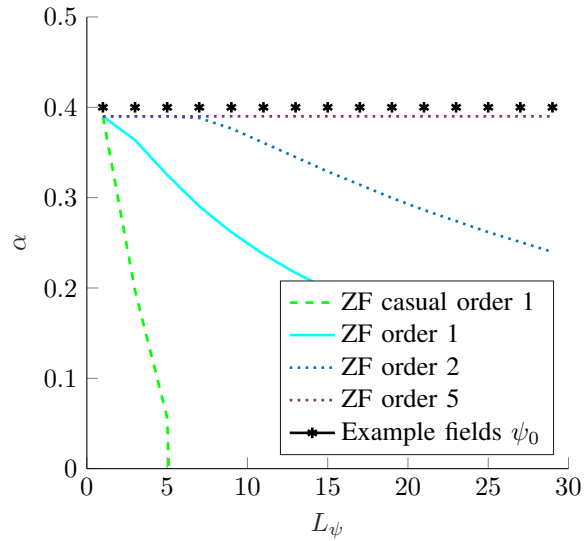

We now consider an LPV system $G(\rho)$ described by,
\begin{equation}\label{eq:lpv_example}
\begin{split}
\dot{x}&=v\\
\dot{v}&=-\rho(t)v -u,
\end{split}
\end{equation}
where, $\rho(t) \in \mathcal{P}=[0.8,1.2] \quad \forall t \in [0,\infty)$.
The scheduling parameter $\rho$ can be seen as a time-varying or adaptive damping co-efficient and can be either fixed and unknown or time-varying. 
It can be verified that \eqref{eq:perf_LMI_lpv} for this example is affine in $\Bar{\rho}$ and hence satisfaction of the inequality for $\Bar{\rho}=0.8$ and $\Bar{\rho}=1.2$ implies the satisfaction for any $\Bar{\rho}\in [0.8,1.2]$ \cite{.2000scherer}.
This reduces the condition \eqref{eq:perf_LMI_lpv} to a finite dimensional feasibility problem that is implemented to produce the results discussed next.
Fig. \ref{fig:LPV_generic_vehicle} shows the convergence rate estimates provided by different multipliers for fields $\psi \in \mathcal{S}(1,L_{\psi})$ with increasing $L_{\psi}$.
As in the previous examples, quadratic fields are chosen to get an upper bound on the convergence rate.
The reduction in conservatism with increasing order of multiplier can be clearly seen.
For this example, fifth-order ZF multipliers show negligible conservatism.
This analysis essentially guarantees, that for this chosen example, poorly conditioned fields do not affect the convergence rate.
For constant trajectories, i.e.,  $\rho(t)=\Bar{\rho} \in [0.8,1.2] \quad \forall t$, and for quadratic fields (linear gradients), a root-locus argument can be used to show that $\alpha=0.4$ for any $L_{\psi}\geq m_{\psi}$.
With the current performance analysis, we can see that this holds even for any non-constant trajectories $\rho$ restricted to the allowable parameter range and for any strongly convex field $\psi \in \mathcal{S}(m_{\psi},L_{\psi})$.
This example also illustrates the benefit of non-causal multipliers over causal ones and the reduction in conservatism with increasing order of the ZF multiplier.

\subsection{Quadrotor with two modes}
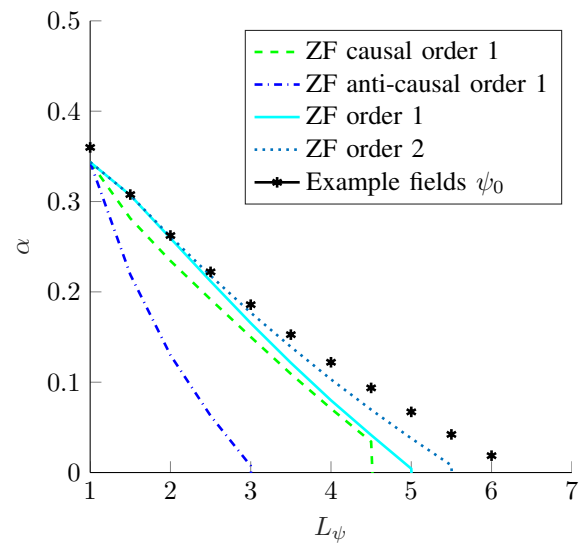
\begin{figure}[t!]
	\centering
%
%
\definecolor{mycolor1}{rgb}{0.00000,1.00000,1.00000}%
\definecolor{mycolor2}{rgb}{0.00000,0.44700,0.74100}%
\definecolor{mycolor3}{rgb}{0.85000,0.32500,0.09800}%
\definecolor{mycolor4}{rgb}{0.92900,0.69400,0.12500}%
\definecolor{mycolor5}{rgb}{0.49400,0.18400,0.55600}%
\begin{tikzpicture}

\begin{axis}[%
width=2.521in,
height=2.366in,
at={(0.758in,0.481in)},
scale only axis,
xmin=1,
xmax=7,
xlabel style={font=\color{white!15!black}},
xlabel={$L_{\psi}$},
ymin=0,
ymax=0.5,
ylabel style={font=\color{white!15!black}},
ylabel={$\alpha$},
axis background/.style={fill=white},
title style={font=\bfseries},
title={},
axis x line*=bottom,
axis y line*=left,
legend style={legend cell align=left, align=left, draw=white!15!black}
]
\addplot [color=green, dashed, line width=1.0pt]
  table[row sep=crcr]{%
1	0.3436279296875\\
1.5	0.2813720703125\\
2	0.234375\\
2.5	0.191650390625\\
3	0.150146484375\\
3.5	0.1092529296875\\
4	0.07080078125\\
4.5	0.0341796875\\
5	-1\\
5.5	-1\\
6	-1\\
6.5	-1\\
7	-1\\
7.5	-1\\
8	-1\\
8.5	-1\\
9	-1\\
9.5	-1\\
10	-1\\
};
\addlegendentry{ZF causal order 1}

\addplot [color=blue, dashdotted, line width=1.0pt]
  table[row sep=crcr]{%
1	0.3436279296875\\
1.5	0.2191162109375\\
2	0.1300048828125\\
2.5	0.0628662109375\\
3	0.00732421875\\
3.5	-1\\
4	-1\\
4.5	-1\\
5	-1\\
5.5	-1\\
6	-1\\
6.5	-1\\
7	-1\\
7.5	-1\\
8	-1\\
8.5	-1\\
9	-1\\
9.5	-1\\
10	-1\\
};
\addlegendentry{ZF anti-causal order 1}

\addplot [color=mycolor1, line width=1.0pt]
  table[row sep=crcr]{%
1	0.3436279296875\\
1.5	0.3070068359375\\
2	0.2593994140625\\
2.5	0.2117919921875\\
3	0.1654052734375\\
3.5	0.1214599609375\\
4	0.0799560546875\\
4.5	0.04150390625\\
5	0.0042724609375\\
5.5	-1\\
6	-1\\
6.5	-1\\
7	-1\\
7.5	-1\\
8	-1\\
8.5	-1\\
9	-1\\
9.5	-1\\
10	-1\\
};
\addlegendentry{ZF order 1}

\addplot [color=mycolor2, dotted, line width=1.0pt]
  table[row sep=crcr]{%
1	0.3436279296875\\
1.5	0.3070068359375\\
2	0.26123046875\\
2.5	0.2178955078125\\
3	0.177001953125\\
3.5	0.13916015625\\
4	0.1031494140625\\
4.5	0.069580078125\\
5	0.037841796875\\
5.5	0.0079345703125\\
6	-1\\
6.5	-1\\
7	-1\\
7.5	-1\\
8	-1\\
8.5	-1\\
9	-1\\
9.5	-1\\
10	-1\\
};
\addlegendentry{ZF order 2}




\addplot [color=black, line width=1.0pt, draw=none, mark=asterisk, mark options={solid, black}]
  table[row sep=crcr]{%
1	0.359877781991861\\
1.5	0.30771931207812\\
2	0.262397425883275\\
2.5	0.222141826976739\\
3	0.185818152359666\\
3.5	0.152649297281869\\
4	0.12207640065513\\
4.5	0.0936823881608809\\
5	0.0671467617957762\\
5.5	0.0422173233107017\\
6	0.0186916779750199\\
6.5	-0.00359531647237843\\
7	-0.0247803316318347\\
7.5	-0.0449778747672073\\
8	-0.0642849197875841\\
8.5	-0.0827843746763419\\
9	-0.100547722380696\\
9.5	-0.117637062795108\\
10	-0.13410671307487\\
};
\addlegendentry{Example fields $\psi_0$}

\end{axis}
\end{tikzpicture}%
	\caption{Convergence rate estimates for quadrotor dynamics with uncertain or switching mass $m\in \{0.2,2\}$ provided by different multipliers (see Remark \ref{rem:CC_ZF_cusal}) for fields $\psi\in \mathcal{S}(1,L_{\psi})$ and convergence rates for $\psi_0(x)=\frac{1}{2}(x-y_{\textnormal{opt}})^T\textnormal{diag}(1,L_{\psi})(x-y_{\textnormal{opt}})$} 
	\label{fig:uncertain_quadrotor}
\end{figure}

We now consider a scenario with a quadrotor, as in Section \ref{eg:quadrotor}, but with two operating modes.
One operating mode corresponds to the quadrotor carrying some load and the other mode corresponds to no-load.
We model this by considering two masses $m \in \{0.2,2\}$ with LQR controllers designed as in Section \ref{eg:quadrotor} for each mode separately.
We consider an arbitrary switching between the two modes and can be modeled as an LPV (or switching) system with $\mathcal{P}=\{1,2\}$ and $\rho(t) \in \mathcal{P} \quad \forall t$.
Fig. \ref{fig:uncertain_quadrotor} shows the convergence rate estimates provided by different multipliers for fields $\psi \in \mathcal{S}(1,L_{\psi})$.
We observe that in comparison to the LTI case (Fig. \ref{fig:quadrotor_robustness} from Section \ref{eg:quadrotor}), the performance is slightly reduced due to the possibility of arbitrary switching between modes.
Furthermore, the estimates with first order ZF multipliers are not tight anymore and we obtain better results with second order ZF multipliers.
No further improvement in the estimates was observed upto 5th order ZF multipliers.

\section{Numerical Examples with Formation Control}\label{sec:numerical_results_formation_control}
We consider again the linearized quadrotor model with an LQR based state-feedback controller tuned for zero steady-state error for step position references. 
As before, the closed-loop system can be represented by \eqref{eq:vehicle_dyn_flock_chap} and it is augmented with dynamics \eqref{eq:vir_vehicle_dyn_flock_chap} to obtain $G$. 
Multiple such quadrotors models are embedded in a field $\psi \in \mathcal{S}(m_{\psi},L_{\psi})$ with an interaction potential $V(y)=\frac{1}{2}(y-r)^T(\mathcal{L} \otimes I_d)(y-r)$, where $\mathcal{L}$ is the graph Laplacian corresponding to a communication graph $\mathcal{G}$ and let $\mathcal{V}_l$ be the set of informed agents. We consider the setup of Problem \ref{prob:formation_control} such that $f$ as defined in \eqref{eq:defn_f_flock} belongs to $\mathcal{S}(m,L)$.
We estimate the rate of convergence for dynamics \eqref{eq:sys_dyn_G_hat} using Theorem \ref{theom:main_perf_analysis_formation} and test the conservatism by finding worst-case examples. 
We first assume that the constants $m$ and $L$ are known in Section \ref{subsec:spec_known} and show how these could be estimated in Section \ref{subsec:spec_unknown}. 
\subsection{Conservatism analysis with known $m$ and $L$} \label{subsec:spec_known}
For a free parameter $L$, consider the uncertainty set $\mathbf{\Delta}_{0.3,L}$ of tuples $(\mathcal{G},\mathcal{V}_l,\psi)$ such that $f$ as defined in \eqref{eq:defn_f_flock} with $V(y)=\frac{1}{2}(y-r)^T(\mathcal{L} \otimes I_d)(y-r)$ belongs to $\mathcal{S}(0.3,L)$.  

\begin{figure}[!t]
	\input{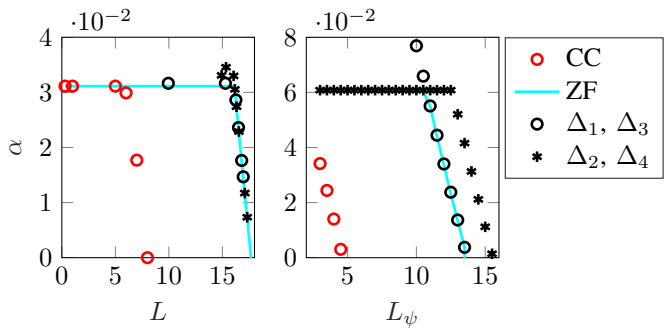}
	\caption{Convergence rate estimates for quadrotor example from Section \ref{sec:numerical_results_formation_control} obtained from Theorem \ref{theom:main_perf_analysis_formation} with circle criterion (CC) IQC and general non-causal Zames-Falb (ZF) IQC. Exact convergence rates for worst-case examples considered in Section \ref{subsec:spec_known} (left) and Section \ref{subsec:spec_unknown} (right) are also shown.}
	\label{fig:Example_spec_known}
\end{figure}

In order to estimate the conservatism, we wish to find non-trivial examples with performance as close as possible to the performance guaranteed by the theory.
For this purpose, define
\begin{equation*}
    \begin{split}
        \Delta_1 & = \{(\mathcal{G},\mathcal{V}_l,\psi):\mathcal{G}=\mathcal{G}_{\textnormal{star}}^5,\mathcal{V}_l=\{1\} , \psi \in \mathcal{S}(m_{\psi},L_{\psi})\},\\
        \Delta_2 & = \{(\mathcal{G}^{25},\mathcal{V}_l,\psi):\psi = 1.85 ||y-y_*||^2 \}.
    \end{split}
\end{equation*}
The graph $\mathcal{G}$ and $\mathcal{V}_l$ is fixed in $\Delta_1$ and constants $m_{\psi}$ and $L_{\psi}$ are chosen such that $\Delta_1 \subseteq \mathbf{\Delta}_{0.3,L}$.
In set $\Delta_2$, $\psi$ is kept fixed while a graph with 25 nodes is chosen such that $\Delta_2 \subseteq \mathbf{\Delta}_{0.3,L}$.
Fig. \ref{fig:Example_spec_known} (left) shows the convergence rate estimates obtained from Theorem \ref{theom:main_perf_analysis_formation} guaranteed by the Zames Falb multipliers and the circle criterion for different values of $L$.
While the circle criterion can certify stability for $0.3\leq L \leq 7$, the Zames-Falb multipliers can certify stability for $0.3\leq L \leq 17.64$.
Since convergence rates for examples in $\Delta_1$ and $\Delta_2$ coincide with the estimates, these estimates are tight.
\subsection{Robust analysis with unknown $m$ and $L$} \label{subsec:spec_unknown}
We now show that constants $m$ and $L$ can be estimated using structural properties of the graph.
We assume that $\psi \in \mathcal{S}(m_{\psi},L_{\psi})$ and a minimal structure in the graph is known in the form of essential edges that are present in all allowable graphs.
This means that all allowable graphs could be obtained by adding edges to the minimal graph.
Let the graph Laplacian associated with this minimal graph be denoted by $\mathcal{L}_0$.
As defined in \ref{defn:Lgs_defn}, let the grounded graph Laplacian associated with this minimal graph be denoted by $\mathcal{L}_m=\mathcal{L}_0+m_{\psi}E$.
Using the fact that adding edges can only increase the eigenvalues of the graph Laplacian (see Lemma 6.9 from \cite{FB-LNS}), we have that 
$\mathcal{L}_m=\mathcal{L}_0+m_{\psi}E\preceq \mathcal{L}+m_{\psi}E= \mathcal{L}_s$.
Furthermore, assuming that the maximum degree, denoted by $d_{\textnormal{max}}$ is known and using the fact that $2d_{\textnormal{max}}$ is an upper bound on the largest eigenvalue of the graph Laplacian (apply Gershgorin Disks Theorem) together with $\psi \in \mathcal{S}(m_{\psi},L_{\psi})$, we get, $\mathcal{L}_b \preceq (2d_{\textnormal{max}}+L_{\psi}) I_N$.
The uncertainty set can then be defined to be $\mathbf{\Delta}_{m,L}$ with $m=\lambda_{\textnormal{min}}(\mathcal{L}_m)$ and $L=(2d_{\textnormal{max}}+L_{\psi})$.

This is now illustrated on a concrete example.
Let $\mathbf{\Delta}$ be the set of all $(\mathcal{G},\mathcal{V}_l,\psi)$ that satisfy the following assumption.
\begin{enumerate}
    \item At least one third of total number of agents are informed agents (have access to the gradient).
	\item Every agent that is not an informed agent has an edge with at least one informed agent.
	\item Maximum degree of all agents is $2$.				
	\item $\psi \in \mathcal{S}(3,L_{\psi})$.
\end{enumerate}
Since any informed agent $i \in \mathcal{V}_l$ is either connected to 0, 1 or 2 other agents, a suitable ordering of the agents will lead to a minimal grounded Laplacian of the form
$\mathcal{L}_m=\textnormal{blkdiag}(m_{\psi}I,\mathcal{L}_1,\cdots,\mathcal{L}_1,\mathcal{L}_2,\cdots,\mathcal{L}_2)$,
where, 
\begin{equation*}
    \begin{split}
        \mathcal{L}_1=\begin{bmatrix}1+m_{\psi}&-1\\-1&1\end{bmatrix},
        \mathcal{L}_2=\begin{bmatrix}2+m_{\psi}&-1&-1\\-1&1&0\\0&-1&1\end{bmatrix}.
        \end{split}
\end{equation*}
Therefore,  we obtain
\begin{equation*}
\begin{split}
    m&=\lambda_{\textnormal{min}}(\mathcal{L}_m)=\textnormal{min}(m_{\psi},\lambda_{\textnormal{min}}(\mathcal{L}_1),\lambda_{\textnormal{min}}(\mathcal{L}_2))=0.4116,\\
    L&=L_{\psi}+2*d_{\textnormal{max}}=L_{\psi}+4.
\end{split}
\end{equation*}
In other words, $\mathbf{\Delta} \subseteq \mathbf{\Delta}_{0.4116,L_{\psi}+4}$. Fig.\ref{fig:Example_spec_known} (right) shows the convergence rate estimates for different values of $L_{\psi}$.
We can again find worst-case examples that coincide with the estimates showing the analysis is without any conservatism here as well.
The worst case examples here turn out to be in the sets $\Delta_3 \subseteq \mathbf{\Delta}$ and  $\Delta_4 \subseteq \mathbf{\Delta}$ defined as
\begin{equation*}
    \begin{split}
        \Delta_3 & = \{(\mathcal{G},\mathcal{V}_l,\psi):\mathcal{G}=\mathcal{G}_{\textnormal{cycle}}^4,\mathcal{V}_l=\mathcal{V} , \psi = \psi_0 \},\\
        \Delta_4 & = \{(\mathcal{G},\mathcal{V}_l,\psi):\mathcal{G}=\mathcal{G}_{\textnormal{star}}^3,\mathcal{V}_l=\{1\}, \psi = \psi_0 \},
    \end{split}
\end{equation*}
where, $\psi_0(x)=(x-y_{\textnormal{opt}})^T\begin{bmatrix}3&\\&L_{\psi}\end{bmatrix}(x-y_{\textnormal{opt}})$.
We emphasize that the robust analysis demonstrated here only requires the knowledge of the maximum degree and a minimal graph.

\section{Numerical Results with Flocking} \label{sec:numerical_results_flocking}
\begin{figure}[t]
	\centering
	\input{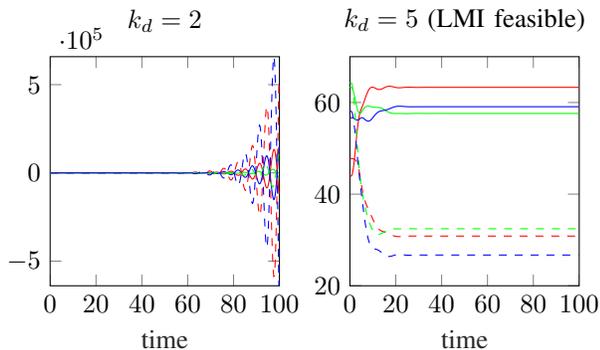}
	\caption{Flocking trajectories (solid lines: x-coordinate and dashed line: y coordinate) for different values of $k_d$ with 3 agents and source located at (60,30)}
	\label{fig:flocking_simulation}
\end{figure}

We continue the quadrotor example from Section \ref{sec:numerical_results_formation_control} and consider three quadrotors embedded in $\psi \in \mathcal{S}(0.5,1)$ with the minimizer $y_{\textnormal{opt}}=\left[60\;\;30\right]^T$.
Let 
$    V(y)=\sum_{(i,j)}k\frac{1}{2}(||y_i - y_j||_{\sigma}-d)^2,$
where, the sigma-norm $||.||_{\sigma}$ is as defined in \cite{OlfatiSaber.2006} (with $\varepsilon=1$) and $k$ is uncertain with $0\leq k \leq 1$.
We assume that at least one agent has access to the gradient.
With these assumptions, one can show that $f$ (see Definition \ref{defn:f_flock}) is proper since $\psi$ is strongly convex, $V_f$  is non-negative and agents that are not informed agents, if any, are connected to the informed agent at all times.
Furthermore, 
$\nabla^2 f(y)=\nabla^2 V_f (y) + \nabla^2 \left(\sum_{i \in \mathcal{V}_l} \psi(y_i)\right)\preceq 2I + I \preceq 3I,$
which implies that $f$ satisfies Assumption \ref{assum:f_flock} with $\mathcal{S}=\mathbb{R}^{Nd}$ and $M_{10}=M_{20}=\begin{bmatrix}9I&0\\0&-I\end{bmatrix}$.
Furthermore, the local LQR controller is designed to satisfy Assumption \ref{assm:ss_error}.
Therefore, applying Corollary \ref{corr:global_flocking_analysis}, a sufficient condition for stability of the overall dynamics is the feasibility of the LMI \eqref{eqn:LMI_flock_Z}.
Considering $k_d$ as the tuning gains, numerical studies verify that the LMI \eqref{eqn:LMI_flock_Z} is feasible for $k_d \geq 4.8$.
Fig. \ref{fig:flocking_simulation} shows trajectories with different values of $k_d$ for $\psi(z)=\frac{1}{2}||z-y_{\textnormal{opt}}||^2$ and $k=1$.
As can be seen from Fig. \ref{fig:flocking_simulation}, the trajectories for $k_d\leq 2$ do not converge, whereas with $k_d=5$, the trajectories converge.
It was found that the trajectories with $k_d=4$ also converge showing a potential conservatism.

\section{Conclusions}
We analyze cooperative control dynamics involving gradient based forcing terms using the framework of IQCs and demonstrate its application on the problem of source-seeking under formation control and flocking dynamics.
Various numerical examples including an LTI quadrotor model with an LQR based controller and two LPV examples where the conservatism can be reduced by increasing the order of the ZF multiplier demonstrate the effectiveness of the approach.
The reduction in conservatism by searching over general higher order non-causal multipliers is also evident from some of the examples considered.
All sufficient conditions developed here are independent of the size of the network owing to the diagonal repeated structure in the multipliers thereby allowing the analysis of large networks.
Flocking dynamics are systematically analyzed for asymptotic stability.

\bibliography{root}

\begin{thebibliography}{10}

\bibitem{Attallah.2020}
Aly Attallah, Adwait Datar, and Herbert Werner.
\newblock Flocking of linear parameter varying agents: Source seeking
  application with underwater vehicles.
\newblock {\em IFAC-PapersOnLine}, 53(2):7305--7311, 2020.

\bibitem{FB-LNS}
Francesco Bullo.
\newblock {\em Lectures on Network Systems}.
\newblock {Kindle Direct Publishing}, 1.5 edition, 2021.

\bibitem{Caverly.20032019}
Ryan~James Caverly and James~Richard Forbes.
\newblock Lmi properties and applications in systems, stability, and control
  theory.

\bibitem{J.Cortes.2017}
Jorge Cort{\'e}s and Magnus Egerstedt.
\newblock Coordinated control of multi-robot systems: A survey.
\newblock {\em SICE Journal of Control, Measurement, and System Integration},
  10(6):495--503, 2017.

\bibitem{datar2023gradient}
Adwait Datar, Antonio~Mendez Gonzalez, and Herbert Werner.
\newblock Gradient-based cooperative control of quasi-linear parameter varying
  vehicles with noisy gradients, accepted for publication at {IFAC} 2023.
\newblock {\em arXiv preprint arXiv:2304.03264}, 2023.

\bibitem{datar_2022_6672230}
Adwait Datar, Christian Hespe, and Herbert Werner.
\newblock {Code for the paper on Robust Performance Analysis of Source-seeking
  Dynamics with Integral Quadratic Constraints: Initial submission}, 2022.

\bibitem{Datar.51220205152020}
Adwait Datar, Peter Paulsen, and Herbert Werner.
\newblock Flocking towards the source: Indoor experiments with quadrotors.
\newblock In {\em 2020 European Control Conference (ECC)}, pages 1638--1643.
  IEEE, 5/12/2020 - 5/15/2020.

\bibitem{datar2021iqc_source_seek}
Adwait Datar and Herbert Werner.
\newblock Robust performance analysis of source-seeking dynamics with integral
  quadratic constraints.
\newblock In {\em 2022 Annual American Control Conference (ACC)}. IEEE, 2022.

\bibitem{Fax.2004}
J.~A. Fax and R.~M. Murray.
\newblock Information flow and cooperative control of vehicle formations.
\newblock {\em IEEE Transactions on Automatic Control}, 49(9):1465--1476, 2004.

\bibitem{Fazlyab.2018}
Mahyar Fazlyab, Alejandro Ribeiro, Manfred Morari, and Victor~M. Preciado.
\newblock Analysis of optimization algorithms via integral quadratic
  constraints: Nonstrongly convex problems.
\newblock {\em SIAM Journal on Optimization}, 28(3):2654--2689, 2018.

\bibitem{Freeman.62720186292018}
Randy~A. Freeman.
\newblock Noncausal zames-falb multipliers for tighter estimates of exponential
  convergence rates.
\newblock In {\em 2018 Annual American Control Conference (ACC)}, pages
  2984--2989. IEEE, 6/27/2018 - 6/29/2018.

\bibitem{hu2017dissipativity}
Bin Hu and Laurent Lessard.
\newblock Dissipativity theory for nesterov’s accelerated method.
\newblock In {\em International Conference on Machine Learning}, pages
  1549--1557. PMLR, 2017.

\bibitem{Hu.2016}
Bin Hu and Peter Seiler.
\newblock Exponential decay rate conditions for uncertain linear systems using
  integral quadratic constraints.
\newblock {\em IEEE Transactions on Automatic Control}, 61(11):3631--3637,
  2016.

\bibitem{UlfJonsson.}
Ulf Jönsson.
\newblock Lecture notes on integral quadratic constraints, 2001.

\bibitem{khong2014multi}
Sei~Zhen Khong, Ying Tan, Chris Manzie, and Dragan Ne{\v{s}}i{\'c}.
\newblock Multi-agent source seeking via discrete-time extremum seeking
  control.
\newblock {\em Automatica}, 50(9):2312--2320, 2014.

\bibitem{Lessard.2016}
Laurent Lessard, Benjamin Recht, and Andrew Packard.
\newblock Analysis and design of optimization algorithms via integral quadratic
  constraints.
\newblock {\em SIAM Journal on Optimization}, 26(1):57--95, 2016.

\bibitem{Magnus.1979}
Jan~R. Magnus and H.~Neudecker.
\newblock The commutation matrix: Some properties and applications.
\newblock {\em The Annals of Statistics}, 7(2):381--394, 1979.

\bibitem{megretski1997system}
Alexandre Megretski and Anders Rantzer.
\newblock System analysis via integral quadratic constraints.
\newblock {\em IEEE Transactions on Automatic Control}, 42(6):819--830, 1997.

\bibitem{michalowsky2016extremum}
Simon Michalowsky and Christian Ebenbauer.
\newblock Extremum control of linear systems based on output feedback.
\newblock In {\em 2016 IEEE 55th Conference on Decision and Control (CDC)},
  pages 2963--2968. IEEE, 2016.

\bibitem{nelson2018integral}
Zachary~E Nelson and Enrique Mallada.
\newblock An integral quadratic constraint framework for real-time steady-state
  optimization of linear time-invariant systems.
\newblock In {\em 2018 Annual American Control Conference (ACC)}, pages
  597--603. IEEE, 2018.

\bibitem{ogren2004cooperative}
Petter Ogren, Edward Fiorelli, and Naomi~Ehrich Leonard.
\newblock Cooperative control of mobile sensor networks: Adaptive gradient
  climbing in a distributed environment.
\newblock {\em IEEE Transactions on Automatic control}, 49(8):1292--1302, 2004.

\bibitem{OlfatiSaber.2006}
R.~Olfati-Saber.
\newblock Flocking for multi-agent dynamic systems: Algorithms and theory.
\newblock {\em IEEE Transactions on Automatic Control}, 51(3):401--420, 2006.

\bibitem{pfifer2015robustness}
Harald Pfifer and Peter Seiler.
\newblock Robustness analysis of linear parameter varying systems using
  integral quadratic constraints.
\newblock {\em International Journal of Robust and Nonlinear Control},
  25(15):2843--2864, 2015.

\bibitem{Scherer.12022021}
Carsten Scherer and Christian Ebenbauer.
\newblock Convex synthesis of accelerated gradient algorithms.
\newblock {\em SIAM Journal on Control and Optimization}, 59(6):4615--4645,
  2021.

\bibitem{scherer_weiland.2000}
Carsten Scherer and Siep Weiland.
\newblock {\em Linear matrix inequalities in control}.
\newblock 2000.

\bibitem{.2000scherer}
Carsten Scherer and Siep Weiland.
\newblock Linear matrix inequalities in control.
\newblock {\em Lecture Notes, Dutch Institute for Systems and Control, Delft,
  The Netherlands}, 3(2), 2000.

\bibitem{Scherer.16052021}
Carsten~W Scherer.
\newblock Dissipativity and integral quadratic constraints: Tailored
  computational robustness tests for complex interconnections.
\newblock {\em IEEE Control Systems Magazine}, 42(3):115--139, 2022.

\bibitem{shamma1992linear}
Jeff~S Shamma and James~R Cloutier.
\newblock A linear parameter varying approach to gain scheduled missile
  autopilot design.
\newblock In {\em 1992 American Control Conference}, pages 1317--1321. IEEE,
  1992.

\bibitem{Shamma.1992}
Jeff~S. Shamma and James~R. Cloutier.
\newblock A linear parameter varying approach to gain scheduled missile
  autopilot design.
\newblock In {\em 1992 American Control Conference}. IEEE, 1992.

\bibitem{sundararajan2017robust}
Akhil Sundararajan, Bin Hu, and Laurent Lessard.
\newblock Robust convergence analysis of distributed optimization algorithms.
\newblock In {\em 2017 55th Annual Allerton Conference on Communication,
  Control, and Computing (Allerton)}, pages 1206--1212. IEEE, 2017.

\bibitem{sundararajan2020analysis}
Akhil Sundararajan, Bryan Van~Scoy, and Laurent Lessard.
\newblock Analysis and design of first-order distributed optimization
  algorithms over time-varying graphs.
\newblock {\em IEEE Transactions on Control of Network Systems},
  7(4):1597--1608, 2020.

\bibitem{veenman2016robust}
Joost Veenman, Carsten~W Scherer, and Hakan K{\"o}ro{\u{g}}lu.
\newblock Robust stability and performance analysis based on integral quadratic
  constraints.
\newblock {\em European Journal of Control}, 31:1--32, 2016.

\bibitem{Veenman.2016}
Joost Veenman, Carsten~W. Scherer, and Hakan K{\"o}ro{\u{g}}lu.
\newblock Robust stability and performance analysis based on integral quadratic
  constraints.
\newblock {\em European Journal of Control}, 31:1--32, 2016.

\bibitem{Xia.2017}
Weiguo Xia and Ming Cao.
\newblock Analysis and applications of spectral properties of grounded
  laplacian matrices for directed networks.
\newblock {\em Automatica}, 80:10--16, 2017.

\bibitem{Zhang.25022019}
Jingfan Zhang, Peter Seiler, and Joaquin Carrasco.
\newblock Noncausal fir zames-falb multiplier search for exponential
  convergence rate.

\end{thebibliography}
\bibliographystyle{plain}

\appendix
\subsection{Supporting Lemma}\label{appendix:Lemma1}
This section proves a lemma that is central in the derivation of the ZF IQCs and is used in the proof of Theorem \ref{theom:theorem_pqw_filtered_ZF}.
It is covered by \cite[Lemma 3]{Freeman.62720186292018} where the result is presented in a very general setting of Bochner spaces.  
Moreover, since the proof of \cite[Lemma 3]{Freeman.62720186292018} is unavailable, we present a self-contained proof here making all arguments in time-domain.
\begin{lemma}\label{theom:lemma_pq_ZF}
	Let $\alpha \geq 0$ be fixed and let $\beta(\tau)=\textnormal{min}\{1,e^{-2\alpha \tau}\}$ for $\tau \in \mathbb{R}$. 
	Let $\Tilde{u},\Tilde{y} \in \mathcal{L}_{2e}[0,\infty)$ be related by $\Tilde{u}=\nabla f(\Tilde{y}+y_*)$, where $f \in \mathcal{S}(m,L)$ and $y_*$ minimizes $f$.
	Then, the signals $p$ and $q$ defined in \eqref{eq:signal_defns} satisfy, $\forall \tau \in \mathbb{R}$, $\forall T \geq 0$, 
	\begin{equation}\label{eq:lemma_ineq}
		\int_0^T e^{2\alpha t} p(t)^T(q(t)-\beta(\tau)q_T(t-\tau))dt \geq 0,
	\end{equation}
	where $q_T$ denotes the extension defined by
	\begin{equation}\label{eq:extension_signal}
		\begin{split}
			q_T(t)&=
			\begin{cases}
				q(t), & \text{if}\ t\in [0,T], \\
				0, & \text{if}\ t\in \mathbb{R} \backslash [0,T].
			\end{cases}
		\end{split}
	\end{equation}
\end{lemma}
\begin{proof}
	The proof goes along the lines of \cite{Lessard.2016} and borrows some ideas from \cite{Scherer.16052021}.
The central idea behind the proof is to exhibit a non-negative function $F:\mathbb{R}\rightarrow\mathbb{R}$ with support $[0,T]$ such that the integrand of \eqref{eq:lemma_ineq} can be lower-bounded for all $\tau \in \mathbb{R}$ and all $t \in [0,T]$ as 
		\begin{equation}\label{eq:temp_lowerbound}
			e^{2\alpha t} p(t)^T(q(t)-\beta(\tau)q_T(t-\tau))\geq F(t) -F(t-\tau). 
		\end{equation}
		If such a function $F$ exists, integrating both sides of \eqref{eq:temp_lowerbound} from $0$ to $T$ and using the non-negativity of $F$ along with the fact that $F$ is zero outside $[0,T]$, we get the desired result \eqref{eq:lemma_ineq}. 
		The remainder of this proof serves to produce this function $F$.

For convenience, let the dimension of $\tilde{y}$ be denoted by $n_y$, i.e., $\tilde{y} \in \mathbb{R}^{n_y}$.
Let $g: \mathbb{R}^{n_y} \xrightarrow[]{} \mathbb{R}$ be defined as
\begin{equation}
	g(\Tilde{y})=f(\Tilde{y}+y_*)-f(y_*)-\frac{m}{2}||\Tilde{y}||^2,
\end{equation}
where $f \in \mathcal{S}(m,L)$ and $y_*$ minimizes $f$.
This can be used to show that $g \in \mathcal{S}(0,L-m)$, $g(\mathbf{0})=0$ and $\nabla g (\mathbf{0})=\mathbf{0}$. 
It can be further shown (see \cite{Lessard.2016}) that for all $\Tilde{y}$, $\Tilde{y}_1$, $\Tilde{y}_2 \in \mathbb{R}^{n_y}$,
\begin{equation}\label{eq:signal_defn_r}
	(L-m)g(\Tilde{y}) -\frac{1}{2}||\nabla g(\Tilde{y})||^2 \geq 0,
\end{equation}
\begin{equation}\label{eq:prop_g}
	(L-m)\nabla g(\Tilde{y})^T\Tilde{y} \geq (L-m)g(\Tilde{y}) +\frac{1}{2}||\nabla g(\Tilde{y})||^2,
\end{equation}
\begin{equation}\label{eq:prop_g2}
	\begin{split}
		\nabla g(\Tilde{y}_1)^T(\Tilde{y}_1-\Tilde{y}_2) &\geq  g(\Tilde{y}_1)-g(\Tilde{y}_2) \\ 
		& \hspace{0.8cm}+ \frac{||\nabla g (\Tilde{y}_1) - \nabla g (\Tilde{y}_2)||^2}{2(L-m)}.
	\end{split}
\end{equation}
Using \eqref{eq:signal_defn_r}, define a non-negative function $r: \mathbb{R}^{n_y} \xrightarrow[]{} \mathbb{R}$ as
\begin{equation}
	r(\Tilde{y})=(L-m)g(\Tilde{y}) -\frac{1}{2}||\nabla g(\Tilde{y})||^2.
\end{equation}
Using definitions \eqref{eq:signal_defns} and $\Tilde{u}=\nabla f(\Tilde{y}+y_*)$, we have 
\begin{equation}
	\begin{split}
		p(t) &= \nabla g(\Tilde{y}(t)),\\
		q(t) &= (L-m)\Tilde{y}(t)-\nabla g(\Tilde{y}(t)).
	\end{split}
\end{equation}
Since $\nabla g (\mathbf{0})=\mathbf{0}$ and $f$ is a static map, we can consider signal extensions $\Tilde{u}_T,\Tilde{y}_T,p_T$ and $q_T$ using \eqref{eq:extension_signal} to have $\forall t \in \mathbb{R}$,
\begin{equation*}\label{eq:signal_defns_estensions}
	\begin{split}
		p_T(t)&=\Tilde{u}_T(t)-m \Tilde{y}_T(t)=\nabla g(\Tilde{y}_T(t)),\\
		q_T(t)&=L \Tilde{y}_T(t) - \Tilde{u}_T(t)=(L-m)\Tilde{y}_T(t)-\nabla g(\Tilde{y}_T(t)).
	\end{split}
\end{equation*}
Therefore, using \eqref{eq:prop_g}, we obtain for an arbitrary $t_1 \in\mathbb{R}$,
\begin{align}\label{eq:pq_lower_bound1}
	&p_T(t_1)^Tq_T(t_1) \\
	&=\nabla g (\Tilde{y}_T(t_1))^T ((L-m)\Tilde{y}_T(t_1)-\nabla g(\Tilde{y}_T(t_1))) \nonumber\\ 
	&=(L-m)\nabla g(\Tilde{y}_T(t_1))^T\Tilde{y}_T(t_1)-||\nabla g(\Tilde{y}_T(t_1))||^2 \nonumber \\
	&\geq(L-m)g(\Tilde{y}_T(t_1)) -\frac{1}{2}||\nabla g(\Tilde{y}_T(t_1))||^2 \nonumber \\
	&=r(\Tilde{y}_T(t_1)).
\end{align}
Along the same lines, using \eqref{eq:prop_g2}, we get for $t_1,t_2 \in \mathbb{R}$,
\begin{align}\label{eq:eq:pq_lower_bound2}
	p_T(t_1)^T[q_T(t_1)-q_T(t_2)] \geq r(\Tilde{y}_T(t_1))-r(\Tilde{y}_T(t_2)).
\end{align}
For any $\beta \in [0,1]$, multiplying \eqref{eq:pq_lower_bound1} by $(1-\beta)$, multiplying \eqref{eq:eq:pq_lower_bound2} by $\beta$ and then adding them together, we get
\begin{equation}\label{eq:pq_prop_lower_bound3}
	p_T(t_1)^T[q_T(t_1)-\beta q_T(t_2)] \geq r(\Tilde{y}_T(t_1))-\beta r(\Tilde{y}_T(t_2)).
\end{equation}
Note that $\beta(\tau)$ as defined in the statement of this Lemma lies in the interval $[0,1]$.
We can therefore use \eqref{eq:pq_prop_lower_bound3} along with non-negativity of $r$ and $\beta(\tau)\leq e^{-2\alpha \tau} \quad \forall \tau \in \mathbb{R}$ to obtain
\begin{align}
	p_T(t)^T[q_T(t)-\beta&(\tau)q_T(t-\tau)] \nonumber \\
	&\geq r(\Tilde{y}_T(t))-\beta(\tau) r(\Tilde{y}_T(t-\tau)) \nonumber\\
	&\geq r(\Tilde{y}_T(t))-e^{-2\alpha \tau} r(\Tilde{y}_T(t-\tau)). 
\end{align}
Multiplying both sides by $e^{2\alpha t}$ and using \eqref{eq:extension_signal}, we get that for all $\tau \in \mathbb{R}$ and $t \in [0,T]$,
\begin{align*}
	e^{2\alpha t} p(t)^T(q(t)-&\beta(\tau)q_T(t-\tau)) \geq F(t)-F(t-\tau)\nonumber \\
\end{align*}
where $F(t):=e^{2\alpha t} r(\Tilde{y}_T(t))$ is the sought function.
\end{proof}
\begin{remark}
			Note that the asymmetry of the function $\beta(s)$ (also present in the IQCs from \cite{Freeman.62720186292018}) arises from the factor $e^{2\alpha t}$ that we wish to have in $\alpha$-IQC formulation. 
			Specifically, the proof of Lemma \ref{theom:lemma_pq_ZF} hinges on taking a convex-combination of two inequalities with co-efficients $\beta(s)$ and $(1-\beta(s))$. Since this requires $\beta(s)\in [0,1]$ for all $s \in \mathbb{R}$, we need to saturate $\beta(s)$ to 1 for negative $s$.			
\end{remark}
\subsection{Model Matrices}\label{appendix:model_matrices}
\begin{align*}
	\mathcal{A}_G&=\begin{bmatrix}
		\hat{A}            & \hat{B}_q & \hat{B}_p \\
		\mathbf{0}   & \mathbf{0} & I_{Nd} \\
		\mathbf{0}  & \mathbf{0} & -k_d I_{Nd} 
	\end{bmatrix},
	&&\mathcal{B}_G=\begin{bmatrix}
		\mathbf{0} \\
		\mathbf{0}\\
		-k_p I_{Nd}
	\end{bmatrix}\\
	\mathcal{C}_{G1}&=\begin{bmatrix}
		\mathbf{0}  & I_{Nd} & \mathbf{0}
	\end{bmatrix},
&&\mathcal{C}_{G2}=\begin{bmatrix}
		-\hat{C}  & I_{Nd} & \mathbf{0}
	\end{bmatrix},\\
\mathcal{A}_0&=\begin{bmatrix}
	{A}            & {B}_q & {B}_p \\
	\mathbf{0}   & \mathbf{0} & I_{d} \\
	\mathbf{0}  & \mathbf{0} & -k_d I_{d} 
\end{bmatrix},
&&\mathcal{B}_0=\begin{bmatrix}
	\mathbf{0} \\
	\mathbf{0}\\
	-k_p I_{d}
\end{bmatrix},\\
\mathcal{C}_{10}&=\begin{bmatrix}
	\mathbf{0}  & I_{d} & \mathbf{0}\\
	\mathbf{0}  & \mathbf{0} & \mathbf{0}
\end{bmatrix},
&&\mathcal{D}_{10}=\begin{bmatrix}
	\mathbf{0}  & \mathbf{0} \\
	I_{d}  & \mathbf{0} 
\end{bmatrix},\\
\mathcal{C}_{20}&=\begin{bmatrix}
	-{C}  & I_{d} & \mathbf{0}\\
	\mathbf{0}  & \mathbf{0} & \mathbf{0}
\end{bmatrix},
&&\mathcal{D}_{20}=\begin{bmatrix}
	\mathbf{0}  & \mathbf{0} \\
	\mathbf{0} & I_{d} 
\end{bmatrix}.
\end{align*}
\subsection{Parameterization of $h$}\label{appendix:h_parameterization}
The parameterization of $h$ in the following discussion proceeds along the lines of \cite{Veenman.2016} (with the same notation).
Observe that $w_1$ and $w_2$ defined in \eqref{eq:signal_defns_w12} can be seen as the outputs of fictitious LTI systems with impulse responses $e^{-2\alpha (t-\tau)}h(t-\tau)$ and $e^{-2\alpha (t-\tau)}h(-(t-\tau))$ and excited by inputs $q$ and $p$, respectively. 
Consider a fictitious LTI systems of order $\nu$ (to be chosen) as follows.
Let 
\begin{equation}
	A_{\nu}=
	\begin{bmatrix} 
		\lambda     & 0         &\dots     & 0 \\
		1           &\lambda    & \ddots   &0  \\
		0           & \ddots    & \ddots   &0  \\
		\vdots      & 0         & 1        &\lambda  
	\end{bmatrix},
	B_{\nu}=
	\begin{bmatrix} 
		1 \\
		0\\
		\vdots\\
		0 
	\end{bmatrix}.
\end{equation}
Choosing the $\lambda$ appropriately is an open-problem as noted in \cite{Veenman.2016} and for all numerical experiments in this paper, $\lambda$ is set to $-1$ and the multiplier order $\nu$ is chosen from the set $\{1,2,\cdots,5\}$ with most applications showing good results already with first order multipliers, i.e., $\nu=1$ (see Section \ref{sec:numerical_results}).
Due to this small size, the numerical conditioning issues typically associated with these Jordon blocks are not observed. 
Let $$Q_{\nu}(t)=e^{A_{\nu}t}B_{\nu}=e^{\lambda t}R_{\nu}^2 \left[ 1 \quad t \quad \hdots \quad t^{\nu-1}\right]^T,$$ where
$$ R_{\nu}=\textnormal{diag}(\frac{1}{\sqrt{0!}},\frac{1}{\sqrt{1!}},\cdots,\frac{1}{(\sqrt{\nu-1)!}}).$$
Parameterize $h$ by variables $P_1\in \mathbb{R}^{1\times \nu}$ and $P_3\in \mathbb{R}^{1 \times \nu}$ as
\begin{equation}\label{eq:h_defn}
	\begin{split}
		h(t)&=
		\begin{cases}
			P_1 Q_{\nu}(-t) & \text{if}\ t<0, \\
			P_3 Q_{\nu}(t) & \text{if}\ t\geq 0.
		\end{cases}
	\end{split}
\end{equation}
Let $\Tilde{\pi}_{\nu}(s)=\left[ 1 \quad \frac{s}{(s-\lambda)^{\nu -1}} \quad \hdots \quad \frac{s^{\nu -1}}{(s-\lambda)^{\nu -1}}\right]^T$ and let $\left[
\begin{array}{c|c}
	\Tilde{A}_{\nu} & \Tilde{B}_{\nu} \\
	\hline
	\Tilde{C}_{\nu} & \Tilde{D}_{\nu}
\end{array}
\right]$ be a state-space realization of $\Tilde{\pi}_{\nu}$.
It has been shown in \cite{Veenman.2016} that if $H$, $P_1$, $P_3$ are such that 
\begin{equation}\label{eq:L1_norm_constraint}
	H+(P_1+P_3)A_{\nu}^{-1}B_{\nu} \geq 0
\end{equation}
and $\exists \mathcal{X}_1, \mathcal{X}_3 \in \mathbb{S}^{\nu-1}$ such that for $i=\{1,3\}$,
\begin{equation}\label{eq:positivity}
	\begin{split}
		(*)
		\begin{bmatrix}
			\mathbf{0}    & \mathcal{X}_i & \mathbf{0} \\
			\mathcal{X}_i & \mathbf{0}    & \mathbf{0} \\
			\mathbf{0} & \mathbf{0} & \textnormal{diag}(P_i)
		\end{bmatrix}
		&
		\begin{bmatrix}
			I                       & \mathbf{0} \\
			\Tilde{A}_{\nu}         & \Tilde{B}_{\nu} \\
			R_{\nu}\Tilde{C}_{\nu}  & R_{\nu}\Tilde{D}_{\nu} \\
		\end{bmatrix}
		\succ0,
	\end{split}
\end{equation}
then, $h$ defined in \eqref{eq:h_defn} satisfies \eqref{eq:zf_impulse_resp_cond}.
Thus, every element in the set 
\begin{equation*}
	\mathbb{P}=\left\{
		\begin{bmatrix}
			\mathbf{0}  & \mathbf{0} & 	H  & -P_3 \\
			\mathbf{0}  & \mathbf{0} &	-P_1^T    & \mathbf{0}\\
			*&*    & \mathbf{0}  & \mathbf{0} \\
			*&*    & \mathbf{0}  & \mathbf{0} 
		\end{bmatrix}|H, P_1, P_3\textnormal{ satisfy \eqref{eq:L1_norm_constraint} \textnormal{ and }\eqref{eq:positivity}}\right\}.
	\end{equation*}
	corresponds to an $h$ satisfying \eqref{eq:zf_impulse_resp_cond}.

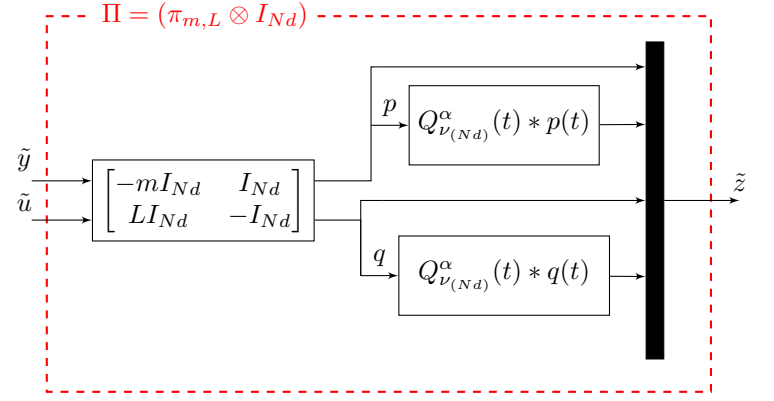
\begin{figure}[t!]
		\centering
		\tikzstyle{block} = [draw, rectangle, 
    minimum height=3em, minimum width=6em]
\tikzstyle{sum} = [draw,  circle, node distance=1cm]
\tikzstyle{input} = [coordinate]
\tikzstyle{output} = [coordinate]
\tikzstyle{pinstyle} = [pin edge={to-,thin,black}]

\begin{tikzpicture}[auto, node distance=2cm,>=latex']
    \node [block] (controller) 
    {$\begin{bmatrix}
    -mI_{Nd}  &  I_{Nd}\\
     LI_{Nd}    & -I_{Nd}
    \end{bmatrix}$};
    \node [block, right of=controller,xshift=2cm,yshift=1cm] (filter1) {$Q_{\nu_{(Nd)}}^{\alpha}(t)*p(t)$};
    \node [block, right of=controller,xshift=2cm,yshift=-1cm,minimum width=2.8cm] (filter2) {$Q_{\nu_{(Nd)}}^{\alpha}(t)*q(t)$};
    \node [input, left=of controller.170,xshift=1.2cm] (input1) {};
    \node [input, left=of controller.190,xshift=1.2cm] (input2) {};
    \node [input, left=of filter1.west,xshift=1.5cm] (p1) {};
    \node [input, left=of filter2.west,xshift=1.5cm] (q1) {};
    \node [output, right of=filter1] (output1) {};
    \node [output, right of=filter2] (output2) {};
    \node [draw,rectangle, right=of controller,xshift=2.4cm,fill=black,minimum height=12em, minimum width=0.05mm] (mux) {};
    \node [input, left=of mux.94.5,xshift=-1.5cm] (p11) {};
    \node [input, left=of mux.180,xshift=-1.5cm] (q11) {};
    \node [output, right=of mux,xshift=-1cm] (z) {};
    \node[text=red, above left= 16mm and -3cm of controller] (Psi) {$\Pi=(\pi_{m,L} \otimes I_{Nd})$};
    \draw[thick,red, dashed] (Psi.east)-|([xshift=-3.8mm]z.west)|-([yshift=-20mm]controller.south)-|([xshift=2mm]input1.west)|-(Psi.west);

    \draw[->] (input1) -- node [name=ytilde,xshift=-0.5cm] {$\Tilde{y}$}(controller.170);
    \draw[->] (input2) -- node [name=utilde,xshift=-0.5cm] {$\Tilde{u}$}(controller.190);
    \draw[-] (controller.10) -| node [name=p] {}(p1);
    \draw[-] (controller.350) -| node [name=q] {}(q1);
    \draw[->] (p1) -- node [name=p] {$p$}(filter1.west);
    \draw[->] (q1) -- node [name=q] {$q$}(filter2.west);
    \draw[-] (p1) |- node [] {}(p11);
    \draw[-] (q1) |- node [] {}(q11);
    \draw[->] (p11) -- node [name=p11] {}(mux.94.5);
    \draw[->] (q11) -- node [name=q11] {}(mux.180);
    \draw[->] (filter1) -- node [name=w2tilde] {}(mux.97.5);
    \draw[->] (filter2) -- node [name=w1tilde] {}(mux.263.5);
    \draw[->] (mux) -- node [name=z,xshift=0.5cm] {$\Tilde{z}$}(z);

\end{tikzpicture}
		\caption{Structure of $\Pi=(\pi_{m,L} \otimes I_{Nd})$.}
		\label{fig:psi_strucure}	
	\end{figure}

	Now, let $A_{\nu}^{\alpha}=A_{\nu}-2\alpha I$, $Q_{\nu}^{\alpha}(t)=e^{-2\alpha t}Q_{\nu}(t)$ and $\pi=\left[\begin{array}{c|c}
		A_{\nu}^{\alpha} & B_{\nu} \\
		\hline 
		0 & 1\\
		I_{\nu} & \mathbf{0}
	\end{array}\right]$.
	With 
\begin{equation*}\label{eq:ss_pi}
	\begin{split}
		\pi_{m,L}=\begin{bmatrix}
			\pi&0\\0&\pi
		\end{bmatrix}\begin{bmatrix}
			-m&1\\L&-1
		\end{bmatrix}
		=
		\left[
		\begin{array}{cc|cc}
			A_{\nu}^{\alpha}         &\mathbf{0}                         &-mB_{\nu}       &B_{\nu} \\
			\mathbf{0}                          &A_{\nu}^{\alpha}        &LB_{\nu}        &-B_{\nu} \\
			\hline
			\mathbf{0}                          &\mathbf{0}                         &-m                      &1\\
			I_{\nu}                  &\mathbf{0}                         &\mathbf{0}                 &\mathbf{0}\\
			\mathbf{0}                          &\mathbf{0}                         &L                       &-1\\
			\mathbf{0}                          &I_{\nu}                 &\mathbf{0}                 &\mathbf{0}
		\end{array}
		\right],
	\end{split}
\end{equation*}
define $\Pi=\pi_{m,L} \otimes I_{Nd}=\left[
\begin{array}{c|c}
	A_{\Pi} & B_{\Pi} \\
	\hline
	C_{\Pi} & D_{\Pi}
\end{array}
\right]$ which is depicted in Fig. \ref{fig:psi_strucure} in the form of a block diagram. 
\subsection{Proofs}\label{appendix:proofs}
\begin{proof}
	[Proof for Lemma \ref{lemm:f_in_Sml}]:
We have, for any ${x},{y} \in \mathbb{R}^{Nd}$,
\begin{equation*}
\begin{split}
(\nabla f({x})-\nabla f({y}))^T ({x} - {y}) &= ({x}-{y})^T(\mathcal{L} \otimes I_d)({x}-{y}) \\
&\hspace{1cm}+({x}-{y})^T{u}_{\psi},\\
\end{split}
\end{equation*}
where $$({x}-{y})^T{u}_{\psi}=\sum_{i \in \mathcal{V}_l} (x_i-y_i)^T(\nabla \psi(x_i)-\nabla \psi(y_i)).$$

Let us first show that statement 2) implies statement 1).
For any $\psi \in \mathcal{S}(m_{\psi},L_{\psi})$, we have that
\begin{equation*}
     m_{\psi} ||x_i-y_i||^2 \leq  (x_i-y_i)^T(\nabla \psi(x_i)-\nabla \psi(y_i)) \leq L_{\psi} ||x_i-y_i||^2.
\end{equation*}
Using Definition \ref{defn:Lgs_defn}, this implies 
\begin{equation} \label{eq: f_ub_Sml}
    (\nabla f({x})-\nabla f({y}))^T ({x} - {y}) \leq ({x}-{y})^T(\mathcal{L}_b \otimes I_d)({x}-{y}),
\end{equation}
\begin{equation} \label{eq: f_lb_Sml}
    (\nabla f({x})-\nabla f({y}))^T ({x} - {y}) \geq ({x}-{y})^T(\mathcal{L}_s \otimes I_d)({x}-{y}).
\end{equation}
Finally, $m I \preceq  \mathcal{L}_s$, $ \mathcal{L}_b \preceq L I$, equations \eqref{eq: f_ub_Sml}, \eqref{eq: f_lb_Sml} and the fact that the spectrum of $X$ and $X \otimes I$ is identical together imply that $f \in \mathcal{S}(m,L)$. 
The reverse direction, i.e., statement 1) implies statement 2), can be shown by choosing the functions
\begin{equation*}
\begin{split}
    \psi_s(z)=m_{\psi}||z||^2 \textnormal{ and }\psi_b(z)&=L_{\psi}||z||^2,
\end{split}
\end{equation*}
defining functions $f_s$ and $f_b$ using Definition \ref{defn:f_flock} for functions $\psi_s$ and $\psi_b$, respectively, 
\begin{equation*}
 \begin{split}
     (\nabla f_s({x})-\nabla f_s({y}))^T ({x} - {y}) &= ({x}-{y})^T(\mathcal{L}_s \otimes I_d)({x}-{y}), \\
     (\nabla f_b({x})-\nabla f_b({y}))^T ({x} - {y}) &= ({x}-{y})^T(\mathcal{L}_b \otimes I_d)({x}-{y}).
 \end{split}
\end{equation*}
This implies statement 2) using the fact that $f_s$, $f_b \in \mathcal{S}(m,L)$.
\end{proof}
\begin{proof}[Proof for Lemma \ref{lemm:m_in_Sml_exists}]:
Let $\mathcal{G}$ and $\mathcal{V}_l$ satisfy Assumption \ref{assm:path_to_leaders}.
This implies that every connected component of $\mathcal{G}$ contains at least one node from $\mathcal{V}_l$.
Therefore, $\mathcal{L}_s$ according to Definition \ref{defn:Lgs_defn}, can be transformed with a permutation of node numbering (similarity transformation), into a block diagonal matrix with each diagonal block being a grounded Laplacian corresponding to each connected component.
Since each diagonal block is a grounded Laplacian with at least one grounded node, the smallest eigenvalue $m_i$ of each diagonal block is positive \cite{Xia.2017}.
By defining $m$ as the minimum over all $m_i$, we obtain $m I \preceq  \mathcal{L}_s$ with $m>0$.
On the other hand, defining $L$ as the maximum eigenvalue of $\mathcal{L}_b$, we obtain $\mathcal{L}_b \preceq L I$.
Using Lemma \ref{lemm:f_in_Sml},  $m I \preceq  \mathcal{L}_s$ and $\mathcal{L}_s \preceq L I $ implies $f \in \mathcal{S}(m,L)$.
This proves the first part of the lemma.

Let us now show the necessity of Assumption \ref{assm:path_to_leaders}.
Without loss of generality, let  $v_1 \in \mathcal{V}$ be such that $\mathcal{G}$ contains no path from $v_1$ to any node in $\mathcal{V}_l$ and let $\mathcal{V}_1=\{v_1,v_2,\cdots,v_{n_1}\}$ be the set of nodes that have a path from $v_1$.
Thus, the grounded Laplacian $\mathcal{L}_s$ as defined in Definition \ref{defn:Lgs_defn} is a block diagonal matrix of the following form
\begin{equation*}
\mathcal{L}_s=
\left[
\begin{array}{cc}
   \mathcal{L}_1 & \mathbf{0} \\
   \mathbf{0}  & \mathcal{L}_{g}
\end{array}
\right],
\end{equation*}
where $\mathcal{L}_1 \in \mathbb{R}^{n_1 \times n_1}$ is a Laplacian without any grounded nodes.
Therefore, $\mathcal{L}\mathbf{1}_{n_1}=\mathbf{0}$.
Hence, a $y \in \mathbb{R}^{N}$ can be constructed as $y=\begin{bmatrix}\mathbf{1}_{n_1}^T & \mathbf{0} \end{bmatrix}^T$, with the property that $||y||^2>0$ and $y^T \mathcal{L}_s y=0$.
Hence, there does not exist an $m>0$ such that $mI \preceq \mathcal{L}_s$.
From Lemma \ref{lemm:f_in_Sml}, this means that if Assumption \ref{assm:path_to_leaders} is violated, there does not exist an $m>0$ such that $f \in \mathcal{S}(m,L)$.
\end{proof}
\begin{proof}[Proof for Lemma \ref{lemm:minimizers}]
First, let ${r}=0$.
For all ${y} \in \mathbb{R}^{Nd}$, we get,
\begin{equation*}
\begin{split}
    f({y})&=\frac{1}{2}{y}^T\mathcal{L}_{(d)}{y} + \sum_{i \in \mathcal{V}_l} \psi(y_i)\geq \sum_{i \in \mathcal{V}_l} \psi(y_i)\geq |\mathcal{V}_l| \psi(y_{\textnormal{opt}}),
\end{split}
\end{equation*}
where the first inequality is obtained by using the fact that the Laplacian $\mathcal{L}$ (and $(\mathcal{L} \otimes I_d)$) is positive semi-definite \cite{FB-LNS} and the second inequality is obtained using Assumption \ref{assum:psi_Sml}.
Since $\mathcal{L}\mathbf{1}_N=\mathbf{0}$, $f(\mathbf{1}_N \otimes y_{\textnormal{opt}})=|\mathcal{V}_l| \psi(y_{\textnormal{opt}})$ which together with the above argument implies $f({y}) \geq f(\mathbf{1}_N \otimes y_{\textnormal{opt}})$ for all ${y} \in \mathbb{R}^{Nd}$, i.e., $\mathbf{1}_N \otimes y_{\textnormal{opt}}$ minimizes $f$.
Since Assumption \ref{assum:psi_Sml} also implies uniqueness of the minimizer, this already completes the proof for part 1).

For part 2), let us assume without loss of generality that $1 \in \mathcal{V}_l$.
Similar to the arguments in 1), for all ${y} \in \mathbb{R}^{Nd}$, we get, 
\begin{equation*}
\begin{split}
    f({y})&=\frac{1}{2}({y}-{r})^T\mathcal{L}_{(d)}({y}-{r}) + \psi(y_1)\geq \psi(y_1)\geq \psi(y_{\textnormal{opt}}).
\end{split}
\end{equation*}
Let $z_1=y_{\textnormal{opt}}$ and $z_j=y_{\textnormal{opt}}+(r_j-r_1)$ as specified in the statement of the lemma.
Thus, ${z}=(\mathbf{1}_N \otimes y_{\textnormal{opt}}) + {r} - (\mathbf{1}_N \otimes r_1)$ and ${z}-{r}=\mathbf{1}_N \otimes (y_{\textnormal{opt}}-r_1)$.
Since $\mathcal{L}\mathbf{1}_N=\mathbf{0}$, $f({z})=\psi(z_1)=\psi(y_{\textnormal{opt}})$ which implies $f({y})\geq f({z})$ for all ${y} \in \mathbb{R}^{Nd}$, i.e., ${z}$ minimizes $f$.
As in part 1), uniqueness of the minimizer (implied by Assumption \ref{assum:psi_Sml}) completes the proof for part 2).

\end{proof}
\begin{proof}[Proof of Lemma \ref{lemm:minimizers_quadratic_radially_symmetric}]
Convexity of $\psi$ implies  that $\nabla \psi(y)(y-y_{\textnormal{opt}})\geq 0$ for all $y \in \mathbb{R}^d$.
This implies that for any $z_i\in \mathbb{R}^d$, $y_{\textnormal{opt}}$ belongs to the half-space $\{y|\nabla \psi(z_i)(z_i-y)\geq 0\}$.
Repeating this argument for each informed agent, we get that $y_{\textnormal{opt}}$ belongs to the intersection of $N$ half-space defined above, i.e., $y_{\textnormal{opt}} \in \{y|\nabla \psi(z_i)(z_i-y)\geq 0 \textnormal{ for all } i \in \mathcal{V}_l\}$.

Assume that $\psi(y)=\psi_r(||y-y_{\textnormal{opt}}||)$, which implies that if $y\neq y_{\textnormal{opt}}$, $\nabla \psi(y)=\frac{\psi_r'(||y-y_{\textnormal{opt}}||)}{||y-y_{\textnormal{opt}}||}(y-y_{\textnormal{opt}})$ and $\nabla \psi(y_{\textnormal{opt}})=0$.
 Additionally, strong convexity of $\psi$ implies $\psi_r'(||y-y_{\textnormal{opt}}||)>0$ for all $y\neq y_{\textnormal{opt}}$.
 Using the fact that $\mathbf{1}_N$ is a left eigenvector of $\mathcal{L}$ with eigenvalue 0, we get
 \begin{equation*}
\begin{split}
    (\mathbf{1}_N^T \otimes I_d)\nabla f ({z})&= \sum_{i \in \mathcal{V}_l}\nabla \psi(z_i) \\
    & = \sum_{i \in \mathcal{V}_l, z_i\neq y_{\textnormal{opt}}}\frac{\psi_r'(||z_i-y_{\textnormal{opt}}||)}{||z_i-y_{\textnormal{opt}}||}(z_i-y_{\textnormal{opt}}).
\end{split}
\end{equation*}
If ${z}$ minimizes $f$, then $\nabla f({z})=0$ and we obtain
\begin{equation}\label{eq:nabla_f1}
\begin{split}
    0& = \sum_{i \in \mathcal{V}_l, z_i\neq y_{\textnormal{opt}}}\frac{\psi_r'(||z_i-y_{\textnormal{opt}}||)}{||z_i-y_{\textnormal{opt}}||}(z_i-y_{\textnormal{opt}})\\
    &=\sum_{i \in \mathcal{V}_l, z_i\neq y_{\textnormal{opt}}}c_i(z_i-y_{\textnormal{opt}})
\end{split}
\end{equation}
where $c_i:=\frac{\psi_r'(||z_i-y_{\textnormal{opt}}||)}{||z_i-y_{\textnormal{opt}}||}$.
Solving for $y_{\textnormal{opt}}$, we get that
\begin{align*}
	y_{\textnormal{opt}}&=\frac{1}{\sum_{i \in \mathcal{V}_l, z_i\neq y_{\textnormal{opt}}} c_i}\sum_{i \in \mathcal{V}_l, z_i\neq y_{\textnormal{opt}}}c_iz_i\\
	&=\sum_{i \in \mathcal{V}_l, z_i\neq y_{\textnormal{opt}}}\mu_iz_i
\end{align*}
where $\mu_i=\frac{c_i}{\sum_{i \in \mathcal{V}_l, z_i\neq y_{\textnormal{opt}}}c_i}$ which implies that $y_{\textnormal{opt}}$ is in the convex hull of $\{z_i|i \in \mathcal{V}_l\}$.

 For part 2) assume that $\psi(y)=\frac{1}{2}y^TQy+b^Ty +c$ which implies $\nabla \psi(y)=Qy + b$ and the unique minimizer is given by $y_{\textnormal{opt}}=-Q^{-1}b$.
Defining $E$ to be a diagonal matrix of size $N$ such that the $i^{th}$ diagonal entry is equal to $1$ if $i\in \mathcal{V}_l$ and equal to $0$ otherwise, we get,
\begin{equation}
    \nabla f({z})=(\mathcal{L} \otimes I_d){z}+(E \otimes Q){z} + (E\mathbf{1}_N) \otimes b.
\end{equation}
Multiplying both sides from the left by $\mathbf{1}_N^T \otimes I_d$, and using the fact that $\mathbf{1}_N$ is a left eigenvector of $\mathcal{L}$ with eigenvalue 0, we get,
\begin{equation*}
\begin{split}
    (\mathbf{1}_N^T \otimes I_d)\nabla f ({z})
    =Q(\sum_{i \in \mathcal{V}_l} z_i)+|\mathcal{V}_l|b.
\end{split}
\end{equation*}
Finally, if ${z}$ is a minimizer of $f$, then $\nabla f({z})=0$ which implies 
$\frac{1}{\mathcal{V}_l}(\sum_{i \in \mathcal{V}_l} z_i)=-Q^{-1}b=y_{\textnormal{opt}}$.

\end{proof}
\begin{proof}[Proof of Theorem \ref{theom:theorem_pqw_filtered_ZF}]
With the signal definitions \eqref{eq:signal_defns_w12}, an appropriate change of integration variable, signal extension as defined in \eqref{eq:extension_signal} and Fubini's theorem, we obtain,
\begin{align}
    &\int_0^T e^{2\alpha t}p(t)^T w_1(t)dt  \\
    &=\int_0^T e^{2\alpha t}p(t)^T \left( \int_{0}^t e^{-2\alpha (t-\tau)}h(t-\tau) q(\tau) d\tau \right)dt \nonumber \\
    &=\int_0^T \int_{s=0}^t e^{2\alpha (t-s)}p(t)^T  h(s) q(t-s) ds dt \nonumber  \\
    &=\int_0^T \int_{s=0}^{\infty} e^{2\alpha (t-s)}p(t)^T  h(s) q_T(t-s) ds dt \nonumber \\
    &=\int_{s=0}^{\infty} h(s) \left(\int_0^T  e^{2\alpha (t-s)}p(t)^T  q_T(t-s)dt \right)  ds.   \label{eq:causal_term} 
\end{align}
Similarly,
\begin{align}
    &\int_0^T e^{2\alpha t}q(t)^T w_2(t)dt\\
    &=\int_{s=-\infty}^0 h(s)  \left(\int_{\tau=0}^T   e^{2\alpha \tau}p(\tau)^Tq_T(\tau-s)  d\tau \right) ds.   \label{eq:anti-causal_term} 
\end{align}
Putting \eqref{eq:causal_term} and \eqref{eq:anti-causal_term} together and using a Lemma \ref{theom:lemma_pq_ZF} (proven in Section \ref{appendix:Lemma1}),
\begin{equation*}
\begin{split}
    &\int_0^T e^{2\alpha t}(p(t)^T w_1(t) + q(t)^T w_2(t)) dt \\
    &=\int_{-\infty}^{\infty} h(s)  \left(  \int_0^T  e^{2\alpha t} \textnormal{min}\{1,e^{-2\alpha s}\} p(t)^T   q_T(t-s)dt \right) ds \\
    &\leq \int_{-\infty}^{\infty} h(s)  \left(  \int_0^T  e^{2\alpha t} p(t)^T   q(t)dt \right) ds \\
    &\leq \int_0^T H e^{2\alpha t} p(t)^T   q(t)dt, \\
\end{split}
\end{equation*}
where the first inequality follows from Lemma \ref{theom:lemma_pq_ZF} and $h(t)\geq 0$ and the second inequality follows from \eqref{eq:zf_impulse_resp_cond}.
\end{proof}
\begin{proof}[Proof of Theorem \ref{theom:theorem_LMI_ZF}]
	From the state-space realization of $\Pi$ (shown in Fig. \ref{fig:psi_strucure}),  the signal definitions \eqref{eq:signal_defns}, \eqref{eq:signal_defns_w12} and \eqref{eq:h_defn} and the block structure of matrices $P \in \mathbb{P}$, we get, 
\begin{equation*}
\begin{split}
 \Tilde{z}^T(t) (P \otimes I) \Tilde{z}(t) =&2(H p(t)^T q(t)-p(t)^T w_1(t)-q(t)^Tw_2(t)),  
\end{split}
\end{equation*}
which can be used together with Theorem \ref{theom:theorem_pqw_filtered_ZF} to finish the proof.

\end{proof}
\begin{proof}[Proof of Theorem \ref{theom:theorem_perf_analysis_lpv}]
Consider the equilibrium $\eta_*,\,u_*=0,\,y_{\textnormal{opt}}$ such that for any trajectory $\rho:[0,\infty) \rightarrow \mathcal{P}$,
\begin{equation} \label{eq:eqm}
	\begin{split}
		0&=A_G(\rho(t))\eta_*, \\
		y_{\textnormal{opt}}&=C_G(\rho(t))\eta_*, \\
		0&=\nabla \psi (y_{\textnormal{opt}}).
	\end{split}
\end{equation}
Note that since the system $G(\rho)$ contains an integrator, following the same arguments as in \cite[Theorem 2.1]{Scherer.12022021}, the existence of an equilibrium $\eta_*$ satisfying the equations \eqref{eq:eqm} is guaranteed by quadratic detectability, i.e., if there exists a matrix $L_o$ such that $(A_G(\bar{\rho})+L_oC_G(\bar{\rho}))$ is Hurwitz for all $\bar{\rho} \in \mathcal{P}$.
Note that the $(1,1)$ block of LMI \eqref{eq:perf_LMI_lpv} ensures that there exists an $\mathcal{X}_G\succ0$ such that
$$A_G(\bar{\rho})^T\mathcal{X}_G+\mathcal{X}_GA_G(\bar{\rho})-C_G(\bar{\rho})^TC_G(\bar{\rho})\prec 0,$$
which implies that the detectability condition \cite[Section 3.12.1]{Caverly.20032019} is verified so that the existence of the equilibrium is established.
Now let $\xi=\begin{bmatrix}
\eta - \eta_* \\
x_{\pi}
\end{bmatrix},$ where $\eta$ is the state of $G(\rho)$, $x_{\pi}$ is the filter state and $\eta_*$ is the equilibrium defined above.
For any trajectory $\rho$, such that, $\rho(t) \in \mathcal{P} \quad \forall t \in [0,\infty)$, the dynamics of $\xi$ can be represented by
\begin{equation} \label{eq:sys_dyn_psi_GI_lpv}
    \begin{split}
        \Dot{\xi}&=\mathcal{A}(\rho(t))\xi + \mathcal{B}(\rho(t))\Tilde{u}, \quad\quad \xi(0)=[\Tilde{\eta}_0^T \quad \mathbf{0}]^T,\\
        \Tilde{z}&=\mathcal{C}(\rho(t))\xi + \mathcal{D}(\rho(t))\Tilde{u}.
    \end{split}
\end{equation}
Since \eqref{eq:sys_dyn_psi_GI_lpv} is the serial interconnection $(\pi_{m,L} \otimes I_d) \begin{bmatrix}
G(\rho)\\
I
\end{bmatrix}$, 
the output $\Tilde{z}$ of \eqref{eq:sys_dyn_psi_GI_lpv} can be obtained from \eqref{eq:signal_defn_z} with signals $\Tilde{u}$ and $\Tilde{y}$, where $\Tilde{y}$ is the output of $G(\rho)$ for input $\Tilde{u}$. 
Furthermore, dynamics  \eqref{eq:sys_dyn_G_lpv} imply that $\Tilde{u},\Tilde{y}$ satisfy $\Tilde{u}=\nabla \psi(\Tilde{y}+y_{\textnormal{opt}})$. 
Hence, Theorem \ref{theom:theorem_LMI_ZF} implies
\begin{equation} \label{eq:zPz_pos}
    \int_0^T e^{2\alpha t}\Tilde{z}^T(t) (P \otimes I_d) \Tilde{z}(t) dt \geq 0 \quad \forall P \in \mathbb{P},\forall T \geq 0.
\end{equation}
Define a storage function $V(\xi)=\xi^T \mathcal{X} \xi$.
Using \eqref{eq:sys_dyn_psi_GI_lpv}, \eqref{eq:perf_LMI_lpv} and the assumption that $\rho(t) \in \mathcal{P} \quad \forall t \in [0,\infty)$, we get
\begin{equation*}
\begin{split}
&\frac{d}{dt}(V(\xi(t)))+2\alpha V(\xi(t))\\
    &=\begin{bmatrix}
    \xi & \Tilde{u}
    \end{bmatrix}
    \begin{bmatrix}
    \mathcal{A}(\rho(t))^T\mathcal{X}+\mathcal{X}\mathcal{A}(\rho(t))+2\alpha \mathcal{X}  & \mathcal{X}\mathcal{B}(\rho(t)) \\
    \mathcal{B}(\rho(t))^T\mathcal{X}    & \mathbf{0}
    \end{bmatrix}
    \begin{bmatrix}
    \xi\\ \Tilde{u}
    \end{bmatrix}\\ 
    &\leq
    -\begin{bmatrix}
    \xi & \Tilde{u}
    \end{bmatrix}
    \begin{bmatrix}
\mathcal{C}(\rho(t))^T \\
\mathcal{D}(\rho(t))^T 
\end{bmatrix}
( P\otimes I_d)
\begin{bmatrix}
\mathcal{C}(\rho(t)) & \mathcal{D}(\rho(t)) 
\end{bmatrix}
\begin{bmatrix}
    \xi\\ \Tilde{u}
    \end{bmatrix}\\
    &=-\Tilde{z}^T(t) (P \otimes I_d) \Tilde{z}(t) .  
\end{split}
\end{equation*}
Rearranging, multiplying by $e^{2\alpha t}$ and integrating from $0$ to $T$, we obtain
\begin{equation*}\label{eq:diff_diss_ineq}
\begin{split}
\frac{d}{dt}(e^{2\alpha t}V(\xi(t))) + e^{2\alpha t}\Tilde{z}^T(t) (P \otimes I_d) \Tilde{z}(t) &\leq 0,\\
e^{2\alpha T}V(\xi(T)) + \int_0^T e^{2\alpha \tau}\Tilde{z}^T(\tau) (P \otimes I_d) \Tilde{z}(\tau) d\tau &\leq V(\xi(0)).
\end{split}
\end{equation*}

Using \eqref{eq:zPz_pos} and $\mathcal{X}succ0$, we get that $V(\xi(t)) \leq e^{-2\alpha T}V(\xi(0))$ implying $||\Tilde{y}(T)|| \leq ||C_G|| \sqrt{\textnormal{cond}(\mathcal{X})} ||\xi(0)||  e^{-\alpha T}$ for all $T\geq0$.
\end{proof}
\begin{proof}[Proof for Lemma \ref{lemm:decomposition}]

It can be shown that there exist permutation matrices $T_1$, $T_2$ of appropriate dimensions such that $T_1^T(X \otimes I)T_2 = I \otimes X$ holds for any real matrix $X$ \cite{Magnus.1979}.
This can be used together with the equivalence of LTI systems under similarity transformations to show that there exist permutation matrices $T$ and $T_{\pi}$ such that, 
\begin{equation}
\begin{split}
    &(\pi \otimes I_{Nd}) \begin{bmatrix}
			I_N \otimes G\\
			I_{Nd}
		\end{bmatrix}
		=(I_{2(\nu+1)} \otimes T_{\pi})(\left[\begin{array}{c|c}
				\mathcal{A}_0     &  \mathcal{B}_0\\
				\hline
				\mathcal{C}_0     &  \mathcal{D}_0
			\end{array}\right]\otimes I_N)T.
\end{split}
\end{equation}
Substituting in LMI \eqref{eq:perf_LMI_big} and using properties of the Kronecker products, the condition reduces to
\begin{equation}
\begin{split}
&\begin{bmatrix}
(\mathcal{A}_0\otimes I_N)^T\mathcal{X}+\mathcal{X}(\mathcal{A}_0\otimes I_N)+2\alpha \mathcal{X}  & \mathcal{X}(\mathcal{B}_0\otimes I_N)T \\
T^T(\mathcal{B}_0\otimes I_N)^T\mathcal{X}    & \mathbf{0}
\end{bmatrix}\\
& \hspace{0.2cm} + (*)
( P\otimes I_d\otimes I_N)
\begin{bmatrix}
(\mathcal{C}_0 \otimes I_N) & (\mathcal{D}_0 \otimes I_N)T 
\end{bmatrix}
\preceq
0, 
\end{split}
\end{equation}
where we have used $(I_{2(\nu+1)} \otimes T_{\pi})^T( P\otimes I_d\otimes I_N)(I_{2(\nu+1)} \otimes T_{\pi})=( P\otimes I_d\otimes T_{\pi}^TT_{\pi})=( P\otimes I_d\otimes I_N)$.
Now, applying a congruence transformation using the permutation matrix $\begin{bmatrix}
I &\mathbf{0}\\
\mathbf{0}&T^T\\
\end{bmatrix}$, we obtain the following equivalent LMI,
\begin{equation} \label{eq:perf_LMI_decomp}
\begin{split}
&\begin{bmatrix}
(\mathcal{A}_0\otimes I_N)^T\mathcal{X}+\mathcal{X}(\mathcal{A}_0\otimes I_N)+2\alpha \mathcal{X}  & \mathcal{X}(\mathcal{B}_0\otimes I_N) \\
(\mathcal{B}_0\otimes I_N)^T\mathcal{X}    & \mathbf{0}
\end{bmatrix}\\
& \hspace{0.2cm} + \left[(*)
( P\otimes I_d)
\begin{bmatrix}
\mathcal{C}_0 & \mathcal{D}_0 
\end{bmatrix} \otimes I_N\right]
\preceq
0. 
\end{split}
\end{equation}
The equivalence between statements 1) and 2) can now be shown by making an argument exactly as in \cite[Section 4.2]{Lessard.2016}.

\end{proof}
\begin{proof}[Proof of Theorem \ref{theom:flocking_analysis_local}]
	With $V_s$ as defined in \eqref{eq: V_flock}, define a function $E:\mathbb{R}\rightarrow \mathbb{R}$ as
\begin{equation*}
	\begin{split}
		E(t) &=V_s({x}(t),{q}(t),{p}(t))\\
		&\geq  ({x(t)}+\hat{A}^{-1}\hat{B_q} {q(t)})^T\hat{Q}({x(t)}+\hat{A}^{-1}\hat{B_q} {q(t)}) \\
		& \hspace{1cm}+ 2\mu (f({q(t)})-f_{\textnormal{min}}).
	\end{split}
\end{equation*}

\noindent Note that if $E(t)\leq V_s({x}_0,{q}_0,{p}_0)$, then  
\begin{equation}
	\begin{split}
		f({q}(t))-f_{\textnormal{min}}& \leq \frac{V_s({x}_0,{q}_0,{p}_0)}{2\mu}  \leq c_1 \textnormal{ and }\\
		||{y}(t)-{q}(t)||^2&=||\hat{C}({x}(t)+\hat{A}^{-1}\hat{B}_q {q}(t))||^2\\
		&\leq||C||^2 \cdot ||{x}(t)+\hat{A}^{-1}\hat{B}_q {q}(t)||^2\\
		&\leq \frac{||C||^2V_s({x}_0,{q}_0,{p}_0)}{\lambda_{\textnormal{min}}(Q)} \leq c_2.
	\end{split}
\end{equation}
This together with the fact that $\mathcal{S}_0$ is bounded and contained in $\mathcal{S}$ (Assumption \ref{assum:f_flock} item 5))  implies that ${q}(t),{y}(t)$ are bounded and contained in $\mathcal{S}$ if $E(t)\leq V_s({x}_0,{q}_0,{p}_0)$.
\normalsize
Now, let $$\left[\begin{array}{c|c}
	\mathcal{A}     &  \mathcal{B}\\
	\hline
	\mathcal{C}_{1}     &  \mathcal{D}_{1} \\
	\mathcal{C}_{2}     &  \mathcal{D}_{2} \\
\end{array}\right]= 
\left[\begin{array}{c|cc}
	\mathcal{A}_G     &  \mathcal{B}_G & -\mathcal{B}_G\\
	\hline
	\mathcal{C}_{G1} & \mathbf{0}     &  \mathbf{0}\\
	\mathbf{0} & I_{Nd} & \mathbf{0}\\
	\mathcal{C}_{G2} & \mathbf{0}     &  \mathbf{0}\\
	\mathbf{0} & \mathbf{0} & I_{Nd}\\
\end{array}\right].$$
Using the structure of the mutipliers $(M_{10} \otimes I_{Nd})$ and $(M_{20} \otimes I_{Nd})$, it can be shown that if there exist $R\succeq 0$, $Q\succ 0$, $\mu>0$, $\lambda_1\geq0$, $\lambda_2\geq0$ such that $\mathcal{Z}$ defined in \eqref{eqn:LMI_flock_Z} is negative semi-definite, then 
$$\mathcal{X}=\begin{bmatrix}
	\hat{\mathcal{X}}_{11}    & \hat{\mathcal{X}}_{12} & \hat{\mathcal{X}}_{13} \\
	\hat{\mathcal{X}}_{21}    & \hat{\mathcal{X}}_{22} & \hat{\mathcal{X}}_{23} \\
	\hat{\mathcal{X}}_{31}    & \hat{\mathcal{X}}_{32} & \hat{\mathcal{X}}_{33} \\
\end{bmatrix}$$
is positive semi-definite and implies
\begin{equation}\label{eqn:LMI_flock_Y}
	\begin{split}
		\mathcal{Y}=\begin{bmatrix}
			{\mathcal{A}}^T{\mathcal{X}}+{\mathcal{X}}{\mathcal{A}}  & {\mathcal{X}}{\mathcal{B}} \\
			{\mathcal{B}}^T{\mathcal{X}}    & \mathbf{0}
		\end{bmatrix}
		+
		\begin{bmatrix}
			\mathbf{0} & \mathbf{0}& \mathbf{0}& \mathbf{0}& \mathbf{0} \\
			\mathbf{0} & \mathbf{0}& \mathbf{0}& \mathbf{0}& \mathbf{0} \\
			\mathbf{0} & \mathbf{0}& \varepsilon I_{Nd}& \mu I_{Nd}&      \mathbf{0} \\
			\mathbf{0} & \mathbf{0}& \mu I_{Nd}&      \mathbf{0}& \mathbf{0} \\
			\mathbf{0} & \mathbf{0}& \mathbf{0}& \mathbf{0}& \mathbf{0} \\
		\end{bmatrix} \\
		+
		\sum_{i=1}^2
		(*)
		(\lambda_i M_{i0} \otimes I_{Nd})
		\begin{bmatrix}
			\mathcal{C}_i & \mathcal{D}_i\\
		\end{bmatrix} \preceq 0.
	\end{split}
\end{equation}
This can be seen by observing that $\mathcal{Y}$ and $\mathcal{Z}$ are block $5 \times 5$ matrices with $\mathcal{Y}_{ij}=I_N \otimes {\mathcal{Z}}_{ij}$ for $i,j \in \{1,\cdots,5\}$. 
Hence, there exists a permutation matrix such that $T^T \mathcal{Y} T = I_N \otimes {\mathcal{Z}}$.
	
Differentiating $E(t)$ with respect to $t$ and using $\mathcal{Y}\preceq0$, we get 
\begin{align*}
	&\dot{E}+  \varepsilon ||{p}||^2+\lambda_1\left((*)\left(M_{10} \otimes I_{Nd}	\right)\begin{bmatrix}
		{q}-q_*\\
		\nabla f({q})
	\end{bmatrix}\right)\\
	&\hspace{1cm}+\lambda_2\left((*)\left(M_{20} \otimes I_{Nd}	\right)\begin{bmatrix}
		{q}-{y} \\
		\nabla f ({q})-\nabla f ({y})
	\end{bmatrix}\right)\leq 0
\end{align*}
for all $t \geq 0$.
Using Assumption \ref{assum:f_flock}, observe that if ${q}(t)$, ${y}(t)$ are in $\mathcal{S}$, then $\dot{E}(t)\leq 0$.
Altogether, if $E(t)\leq V({q}_0,{q}_0,{p}_0)$, then ${q}(t),{y}(t) \in \mathcal{S}_0\subseteq \mathcal{S}$ which implies $\dot{E}(t)\leq 0$ and finally $E(t+s)\leq E(t) \leq V({x}_0,{q}_0,{p}_0)$ for all $s\geq 0$. 
In other words, $\mathcal{S}_0$ is an invariant set.
Finally, observe that ${q}_0,{y}_0$ are in $\mathcal{S}_0$ and hence $\dot{E}(t)\leq -\varepsilon ||{p}(t)||^2$ for all $t\geq 0$.
Applying LaSalle's invariance principle, the trajectory $\bar{\eta}(t)$ converges to the set $\{\begin{bmatrix}
	(-\hat{A}^{-1}\hat{B}_q y_*)^T& y_*^T & 0
\end{bmatrix}^T|\nabla f(y_*)=0,\,y_* \in \mathcal{S} \}$ and since $-\hat{C}\hat{A}^{-1}\hat{B}_q=I$, ${y}(t)$ converges to the set $\{ y_*\in \mathcal{S} |\nabla f(y_*)=0\}$. 
\end{proof}
\end{document}